\documentclass[10pt]{amsart}
\usepackage[centertags]{amsmath}
\usepackage{amsfonts}
\usepackage{amssymb}
\usepackage{amsthm}
\usepackage{graphicx}
[1994/12/01]% LaTeX date must December 1994 or later
\title{$w$-function of the KdV hierarchy}

\thanks{The work was supported  by Russian Foundation for Basic
Research, grant 02-01-0659 and grant for the Leading Russian
Scientific School 2185.2003.1 }

\author{V.M. Buchstaber, S.Yu. Shorina }

\copyrightinfo{2004}%            % copyright year
{American Mathematical Society}% copyright holder
\def\oa#1{\gdef\@oa{#1}}
\def\@setoa{\@oa\@addpunct.}
\def\@maketitle{%
  \normalfont\normalsize
  \let\@makefnmark\relax  \let\@thefnmark\relax
  \ifx\@empty\@oa\else \@footnotetext{\@setoa}\fi
  \ifx\@empty\@subjclass\else \@footnotetext{\@setsubjclass}\fi
  \ifx\@empty\@keywords\else \@footnotetext{\@setkeywords}\fi
  \ifx\@empty\thankses\else \@footnotetext{%
    \def\par{\let\par\@par}\@setthanks}\fi
  \@mkboth{\@nx\shortauthors}{\@nx\shorttitle}%
  \global\topskip8pc\relax
  \@settitle
  \ifx\@empty\authors \else \@setauthors \fi
  \ifx\@empty\@dedicatory
  \else
    \baselineskip18\p@
    \baselineskip26\p@
    \vtop{\centering{\footnotesize\itshape\@dedicatory\@@par}%
      \global\dimen@i\prevdepth}\prevdepth\dimen@i
  \fi
  \@setabstract
  \normalsize
  \if@titlepage
    \newpage
  \else
    \dimen@34\p@ \advance\dimen@-\baselineskip
    \vskip\dimen@\relax
  \fi
} % end \@maketitle
\newcommand{\lcandify}{%
  \nxandlist{\unskip, }{\unskip{} \lowercase{and}~}{\unskip,
\lowercase{and}~}}
\def\@settranslators{\par\begingroup
  \addvspace{6\p@\@plus9\p@}%
  \hbox to\columnwidth{\hss\normalfont\Small
    Translated by %
    \lcandify\@translators \uppercasenonmath\@translators
    \@translators}
  \endgroup
}

\def\enddoc@text{%
  \ifx\@empty\addresses \else\@setaddresses\fi
  \ifx\@empty\@translators \else\@settranslators\fi
}
\newtheorem{thm}{Theorem}[section]
\newtheorem{cor}{Corollary}[section]
\newtheorem{lem}{Lemma}[section]
\newtheorem{ex}{Example}
\newtheorem{defn}{Definition}[section]
\theoremstyle{remark}
\newtheorem{rem}{Remark}[section]
\newtheorem{problem}{Problem}
\numberwithin{equation}{section}

\usepackage[cp1251]{inputenc}  %% 1
\newcommand{\Real}{\mathbb R}

\newcommand{\dd}{\mathrm{ d}}

\newcommand{\prf}{\noindent {\bf Proof.}\;\;}
\newcommand{\ddt}{\partial/\partial t}

\begin{document}

\maketitle

 \begin{abstract}
In this paper we construct a family of commuting multidimensional
differential operators of order $3$, which is closely related to
the KdV hierarchy. We find a common eigenfunction of this family
and an algebraic relation between these operators. Using these
operators we associate a hyperelliptic curve to any solution of
the stationary KdV equation. A basic generating function of the
solutions of stationary KdV equation is introduced as a special
polarization of the equation of the hyperelliptic curve. We also
define and discuss the notion of a $w$-function of a solution of
the stationary $g$-KdV equation.
 \end{abstract}

\section*{Introduction}
At the present time various forms of solutions of the stationary
$g$-KdV equations are known, including the representations with
the $\tau$-function (\cite{hirota}), $\theta$-function
(\cite{ItsMatveev, krich8}), and $\sigma$-function (\cite{1,
Gordon}); rational solutions can be expressed in terms of
Adler-Moser polynomials (\cite{AdlerMoser}). All these functions
satisfy the equation
 %*
 \begin{equation}\label{tau1}
2 \partial_x^2 \log f = - u,
 \end{equation}
 %*
where $u = u(x,t_2, \ldots,t_g)$ is a solution of the stationary
$g$-KdV equation.

In this paper we construct a family of commuting multidimensional
differential operators of third order starting with an arbitrary
solution of the stationary $g$-KdV equation. Using these operators
we solve the following well-known

 \begin{problem}
Supplement \eqref{tau1} with natural conditions so to obtain a
problem with the unique solution.
 \end{problem}

We call this solution a $w$-function of the KdV hierarchy.

In \cite{novikov1} S.P.\,Novikov observed  that each solution of
the stationary $g$-KdV equation is a $g$-gap potential of the
Schr\"{o}dinger operator. It was shown in \cite{1}, \cite{Gordon}
that the Kleinian $\sigma$-function $\sigma(x,t_2,...,t_g)$
provides a solution of the $g$-KdV equation. This fact follows
from a general result describing all algebraic relations between
the higher logarithmic derivatives of the $\sigma$-function.

We are going to discuss also the following natural

 \begin{problem}
Describe all the relations  between the higher logarithmic
derivatives
 %*
 \begin{equation*}
\frac{\partial^{i_1 + \ldots + i_g}}{\partial x^{i_1} \partial
t_2^{i_2} \ldots \partial t_g^{i_g}} \log w(x ,t_2, \cdots,
t_g),\quad \text{where}\; i_1+\ldots+i_g \ge 2
 \end{equation*}
 %*
following from the construction of the $w$-function of the KdV
hierarchy.
 \end{problem}
A solution of this problem is given in Section \ref{basic}.

In the paper \cite{krich} I.M.\,Krichever introduced a concept of
the Baker--Akhiezer function as a common eigenfunction of the
operators $\mathcal{L}$ and $A$ (see section~1 for definitions).
This function is characterized by its analytic properties,
including the behavior at the singular points. In Section
\ref{homcond} we express this function via the common
eigenfunction of our family of commuting differential operators.

The  results of this paper were  partially announced in
\cite{BS1}, \cite{BS2}.

\section{Preliminaries}

This section is a brief review of basic facts about the KdV
hierarchy. See \cite{ZMNP} for more details.

The classical  KdV (Korteweg -- de Vries) equation is
 %*
 \begin{equation}\label{kdv}
\frac{\partial}{\partial t}\,u=\frac{1}{4}(u'''-6 u u'),
 \end{equation}
 %*
where $u$ is a function of real variables $x$ and $t$; the prime
means derivation with respect to $x$.

Denote $\mathcal{L}=\partial_x^2 - u$ the Schr\"{o}dinger operator
with the potential $u$. The second term here means {\em the
operator of multiplication} by the function $u$; we will use
similar notation throughout the paper. Let also
 %*
 \begin{equation}\label{A1}
A_1 = \partial_x^3 - \frac34 (u \partial_x + \partial_x u) =
\partial_x^3 - \frac32 u \partial_x - \frac34 u'.
 \end{equation}
 %*
Then, as it was first noticed in \cite{lax}, the KdV equation is
equivalent to the condition
 %*
 \begin{equation*}
[A_1,\mathcal{L}] = \frac14(u'''-6 u u').
 \end{equation*}
 %*

Denote $\mathfrak D$ a ring of differential operators with
coefficients in the ring of smooth functions in variables $x$ and
$t$. Consider an action of the operator $\ddt$ on the ring
$\mathfrak D$ defined by the formula
 %*
 \begin{equation}\label{dB}
\frac{\partial}{\partial t} \left(\sum_{k \ge 0} f_k(t,x)
\partial_x^k \right) =  \sum_{k \ge 0} \frac{\partial
f_k(t,x)}{\partial t}\,
\partial_x^k.
 \end{equation}
 %*
Then for the operator  $\mathcal{L}$ we obtain the equality
 %*
 \begin{equation*}
\frac{\partial}{\partial t} \mathcal{L} =
-\frac{\partial}{\partial t}\,u.
 \end{equation*}
 %*

So, equation \eqref{kdv} is equivalent to
 %*
 \begin{equation*}
\frac{\partial}{\partial t} \mathcal{L} =[A_1,\mathcal{L}].
 \end{equation*}
 %*

For every differential operator $B \in \mathfrak D$ define its
formal conjugate $B^*$ as follows: take, by definition,
 %*
 \begin{equation}
\partial_x^*=-\partial_x,\quad f^*=f,
 \end{equation}
 %*
where $f$ is an operator of multiplication by the function $f$,
and assume $*$ to be a ring anti-homomorphism:
 %*
 \begin{equation}
 \label{adj}
(B_1 B_2)^* = B_2^* B_1^*,\quad (B_1+B_2)^*=B_1^*+B_2^*,
 \end{equation}
 %*
for all $B_1, B_2 \in \mathfrak D$.

We call an operator  $B$ {\em symmetric} if $B^*=B$, and {\em
anti-symmetric} if $B^*=-B$. Thus, the operator $\mathcal{L}$ is
symmetric. while the operator $A_1$ is anti-symmetric.

Consider a subring $\mathfrak D_1 \subset \mathfrak D$ generated
by $\partial_x$ and the multiplication operator $u$. Supply the
ring $\mathfrak D_1$ with a grading such that
%*
 \begin{equation}\label{deg}
 \deg u=2,\quad \deg \partial_x = 1.
 \end{equation}
%*
Thus, $\deg  u^{(k)} = \deg \partial_x^k u = k+2$. The operators
$\mathcal{L}$ and $A_1$ are then homogeneous of order $2$ and $3$,
respectively.

 \begin{defn}
Denote by $\mathfrak{A}$ the linear space of anti-symmetric
differential operators $A$ such that the commutator
$[A,\mathcal{L}]$ is an operator of multiplication by a function.
 \end{defn}

 \begin{thm}\cite{lax}\label{thmAbasis}
The space $\mathfrak{A}$ has a basis $A_0, A_1, \dots$ where
$A_k=\partial_x^{2k+1}+\sum P_{k,i} \partial_x^i$ is a homogeneous
differential operator of order $2k+1$, and $P_{k,i}$ is a
differential polynomial of $u$ of order $2k+1-i$.
 \end{thm}

The recurrence relation for the operators $A_k$ can be found in
\cite{GelfandDikij}. The operator $A_1$ is given by \eqref{A1}.
The operators $A_0$ and $A_2$ are
 %*
 \begin{align*}
&A_0= \partial_x\\
&A_2 = \partial_x^5 - \frac{5}{4}(u \partial_x^3 + \partial_x^3 u
) + \frac{15}{8} u \partial_x u + \frac{5}{16}(u'' \partial_x +
\partial_x u'').
 \end{align*}
 %*

Denote $r_k[u] = [A_k, \mathcal{L}]$, so that $r_1[u] =
\frac14(u'''-6 u u')$, $r_2[u] = \frac1{16}(u^{(5)}- 10 u u''' -
20 u' u''+ 30 u^2 u')$, etc. Suppose now that $u$ depends on $x$
and an {\em infinite} set of variables $t_1, t_2, \dots$. The
equation
 %*
 \begin{equation}\label{hkdv}
\partial_{t_g} u = r_g[u]
 \end{equation}
 %*
is called the $g$-th higher KdV equation.

The family of equations \eqref{hkdv} is called the KdV hierarchy.

The action of differential operators $\partial_{t_k}$ on the ring
$\mathfrak D_u$ is defined similar to \eqref{dB}.

 \begin{lem}
The operators $A_k$ satisfy the following ``zero curvature''
condition: $\partial_{t_k} A_m-\partial_{t_m} A_k=[A_k,A_m]$, or,
equivalently, $[\partial_{t_k}-A_k,\partial_{t_m}-A_m]=0$.
 \end{lem}

The expression
 %*
 \begin{equation}\label{Lenard}
\mathcal{R}=\frac{1}{4} \partial_x^2 - \frac{1}{2} u'
\partial_x^{-1} -u,
 \end{equation}
 %*
is called the Lenard operator; here $\partial_x^{-1}$ is an
operator of integration with respect to $x$. Note that the Lenard
operator $\mathcal{R}$ is multi-valued, to fix its value we need
to choose the integration constants.

\begin{thm}
Functions $r_k[u]$ are related by the Lenard operator:
 %*
 \begin{equation*}
r_{k+1}[u] = \mathcal{R}(r_k[u]).
 \end{equation*}
 %*
\end{thm}
For example, $r_0=u'$, and $r_1 = \frac14 u''' - \frac32 u u' =
\frac14
\partial_x^2 u' - \frac12 u' u - u u' = \mathcal{R} (r_0)$.

\begin{defn}
The equations
 %*
\begin{equation}\label{novikov}
r_g[u]+\sum_{k=0}^{g-1} a_k r_k[u]=0,
 \end{equation}
 %*
where $a_k$ are constants, are called higher stationary $g$-KdV
(or Novikov) equations.
\end{defn}

 \begin{thm}
A function $u$ is a solution of \eqref{novikov} if and only if it
satisfies the relation $\mathcal{R}^g(u') = 0$ for some choice of
the integration constants; this choice depends on the constants
$a_k$.
 \end{thm}
See \cite{GelfandDikii2} for proof.

If one replaces the function $u$ with $u+c$ where $c$ is a
constant, then the operator $A_k$ becomes $A_k + c
\sum_{i=0}^{k-1} c_{k;i} A_i$ for some constants $c_{k;i}$ where
$c_{k;k-1} \ne 0$. We can choose the constant $c$ so that to
achieve the equality $a_{g-1}=0$ in the decomposition of the
operator $A = A_g+\sum_{k=0}^{g-1} a_k A_k$.

\section{A family of commuting multidimensional differential operators of
order $3$}

In this section we describe a family $\{\mathcal{U}_k\}$ of
differential operators commuting with each other and with  the
Schr\"{o}dinger operator. In this aspect they resemble operators
$\partial_{k+1}-A_k$, but unlike $\{A_k\}$ they are
multidimensional operators of the third order.

Let $\{u_1,u_2,\dots,u_g\}$ be a sequence of functions of
variables $t_1=x,t_2,\dots,t_g$. Denote
$\partial_i=\partial/\partial t_i$. Suppose that the first
derivatives of the function $u_1$ are linearly independent, i.e.\
$\sum_{i=1}^g c_i \partial_i (u_1) \equiv 0$ only if $c_1 = \dots
= c_g = 0$. This condition means that the function
$u_1=u_1(t_1,\dots,t_g)$ essentially  depends on all its
arguments, i.e.\ there is no linear projection  $\pi:\mathbb{C}^g
\to \mathbb{C}^{g-1}$ such that  $u = \pi^* \widetilde{u}$, where
$\widetilde{u}$ is a function on $\mathbb{C}^{g-1}$.

Denote
 %*
 \begin{align*}
&\mathcal{L} = \partial_x^2-u_1,\\
&\mathcal{A}_k = \partial_x^2\partial_k - \frac{1}{2}(u_1
\partial_k +
\partial_k u_1)- \frac{1}{4} (u_k \partial_x +  \partial_x u_k)
 \end{align*}
 %*
where $k = 1, \dots, g$. Note that the $\mathcal{A}_1$ coincides
with the operator $A_1$ given by \eqref{A1}.

Define the formal conjugation $^*$ on the space of
multidimensional differential operators by formulas \eqref{adj}
together with the rule $\partial_i^*=-\partial_i$ where $i>1$. The
operators $\mathcal{A}_k$ are anti-symmetric: $\mathcal{A}_k^*=-
\mathcal{A}_k$.

Consider in the ring of differential operators in variables $t_1,
\dots, t_g$ a subring $\mathfrak D_g$ generated by the operators
$\partial_1, \dots,
\partial_g$ and $u_1, \dots, u_g$. Define on $\mathfrak D_g$ a grading using
formulas \eqref{deg} and assuming also that $\deg u_k=2k$ and
$\deg\partial_k = 2k-1$. It is clear that the operators
$\mathcal{L}$ and $\mathcal{A}_k$ are homogeneous, $\deg
\mathcal{L}=2$, $\deg \mathcal{A}_k = 2k+1$.

 \begin{lem}
The commutator $[\mathcal{L},\mathcal{A}_k]$ is a multiplication
operator if and only if $u'_k=\partial_k u_1$ for all $k$. If this
condition is satisfied then
 %*
 \begin{equation}\label{LAk}
\left[\mathcal{L}, \mathcal{A}_k\right] = \frac14 (u_k'''-2 u_1'
u_k -4 u_1 u_k').
 \end{equation}
 %*
 \end{lem}

\prf Lemma follows from the formula
 %*
 \begin{equation*}
[\mathcal{L},\mathcal{A}_{k}] = (-2 u_k'' + 2 \partial_k u_1')
\partial_x + (- u_k'+ \partial_k u_1)\partial_x^2+\partial_k u_1''
-u_1 \partial_k u_1 - \frac{1}{2} u_k u_1'-\frac{3}{4} u'''_k.
\qed
 \end{equation*}
 %*

 \begin{lem}
Let $K$ be an anti-symmetric differential multidimensional
operator of order $3$. Suppose the commutator $[K, \mathcal{L}]$
is a multiplication operator, and the coefficients in derivatives
of order $3$ are constants. Then $K = \sum_{1\le i\le g} c_i
\mathcal{A}_i + \psi_i \partial_i + \phi$, where $c_i$ are
constants, and the functions $\psi_i$ and $\phi$ do not depend on
$x$.
 \end{lem}

\prf Let $K=\sum_{1 \le i,j,k \le g} s_{ijk}\, \partial_i
\partial_j
\partial_k + \sum_{1 \le i,j \le g} f_{ij}\,\partial_i \partial_j + \sum_{1
\le i \le g} g_i\,\partial_i + h$, where $s_{ijk}$ are constants
such that $s_{ijk} = s_{ikj} = s_{jik}$, and all the $f_{ij}$,
$g_i$, $h$ are functions of $t_1, \dots, t_g$. The anti-symmetry
implies that $f_{ij}=0$ for all $i,j$, and $\sum_{1 \le i \le g}
\frac{\partial g_i}{\partial t_i} = 2 h$.

We have
 %*
 \begin{multline*}
\left[\mathcal{L}, K \right] = \sum_{1\le i \le g} g''_i
\partial_i + 2 \sum_{1\le i \le g} g'_i \partial_x \partial_i +
h'' + 2 h'\partial_x +
\sum_{1 \le i \le g } g_i \frac { \partial}{ \partial t_i } u_1 + \\
\sum_{1 \le i,j,k \le g} s_{ i j k}\left( \frac { \partial^3}{
\partial t_i
\partial t_j
\partial t_m } u_1 +
3 \frac { \partial^2}{ \partial t_i \partial t_j } u_1  \partial_k
+ 3 \frac { \partial}{ \partial t_i  } u_1  \partial_j \partial_k
\right)
 \end{multline*}
 %*

Since the commutator  $\left[\mathcal{L}, K \right]$ is a
multiplication operator, the coefficients at $\partial_i
\partial_j$ and $\partial_i$ in the last formula are zeros.

It follows from the linear independence for the first derivatives
of the function $u_1$ that $s_{ijk}=0$ when $i,j\ne 1$. If $i\ne
1$ then one has $2 g'_i = -3 s_{11i} \frac{\partial}{\partial x}
u_1$.

Equalizing the coefficient at $\partial_x$ to zero, we obtain
$g_1''+2 h'+3 \sum_{1 \le i \le g} s_{1 1 i} \frac{\partial^2
u}{\partial_i \partial x} =0$. Put $c_i=3
s_{11i}=s_{11i}+s_{1i1}+s_{i11}$; $\psi_1=g_1+1/2\,\sum_{i=1}^g
c_i u_i$,  $\psi_k=g_k+ \sum_{i=1}^g c_i u_1$ for $k \ne 1$, and
$\phi=h+\sum_{i=1}^g c_i u_i'$. Then the functions $\psi_i$,
$i=1,\dots,g$ and $\phi$ do not depend on $x$ and $K=\sum_{1\le
i\le g} c_i \mathcal{A}_i+ \psi_i \partial_i+\phi$.\qed

Denote
 %*
 \begin{align*}
&\mathcal{U}_i = \mathcal{A}_i - \partial_{i+1}, \quad\text{for } i<g;\\
&\mathcal{U}_g = \mathcal{A}_g.
 \end{align*}
 %*
The operators $\mathcal{U}_i$ are anti-symmetric and homogeneous.

 \begin{lem}\label{lem2}
The following conditions are  equivalent
 \begin{enumerate}
\item $[\mathcal{L},\mathcal{U}_k]=0$.

\item $-\partial_{k+1} \mathcal{L} = [\mathcal{L},\mathcal{A}_k]$.

\item\label{4qi+1} $\partial_{k+1} u_1 = u_{k+1}' =  \frac14
(u_k''' - 2u_1'u_k - 4u_1 u_k')$ for $k < g$, and $(u_g''' - 2u_1'
u_{g} - 4 u_1 u_g') = 0$.
 \end{enumerate}
 \end{lem}

\prf Lemma follows from the equality
 %*
 \begin{equation*}
[\mathcal{L},\mathcal{U}_k] = [\mathcal{L},\mathcal{A}_k] +
\partial_{k+1} \mathcal{L} = [\mathcal{L},\mathcal{A}_k] -
\partial_{k+1} u_1.\qed
 \end{equation*}
 %*

The last statement in Lemma \ref{lem2} allows to express the
functions $u_k$ recursively via $u_1$ and its $x$-derivative, up
to the choice of a function that do not depend on $x$.

 \begin{cor}
Under the hypotheses of Lemma \ref{lem2}, functions $u_i$ are
related by the Lenard operator $\mathcal{R}$ (see \eqref{Lenard}):
 %*
 \begin{equation}\label{4qi}
\partial_{i+1} u_1 = \mathcal{R}(u_i') = \mathcal{R}^{i+1}(u_1').
 \end{equation}
 %*
\end{cor}

For $i=g$ one has $0 = \partial_{g+1} u_1 =
\mathcal{R}^{g}(u_1')$.

 \begin{cor}
The operators $\{\mathcal{U}_k\}$ commute with $\mathcal{L}$ if
and only if the function $u_1(x)$ is a solution of the stationary
$g$-KdV equation.
 \end{cor}

\begin{lem} \label{lem3}
The operators $\mathcal U_i, \mathcal U_j$ commute for all $1\le
i,j, \le g$ if and only if the functions $\{u_i\}$ satisfy
condition \eqref{4qi+1} of Lemma \ref{lem2} and the following
equalities:
 %*
 \begin{eqnarray}
&\partial_i u_j = \partial_j u_i, \label{diqj-djqi}\\
&u_k' u_i -  u_i' u_k + 2 \partial_{i+1} u_{k}- 2\partial_{k+1}
u_{i} = 0, \quad 1 \le i,k \le g. \label{di+1}
 \end{eqnarray}
 %*
\end{lem}
Lemma is proved by direct calculation.

Note that \eqref{diqj-djqi} implies the existence of a function
$z(t_1,\ldots,t_g)$ such that $\partial_i z = u_i$.

The commutation of operators $\mathcal{U}_i$ is equivalent to zero
curvature conditions for the operators $\mathcal{A}_i$:
 %*
 \begin{equation}\label{nullcur}
\partial_{j+1} \mathcal{A}_i -  \partial_{i+1} \mathcal{A}_j +
\left[\mathcal{A}_i, \mathcal{A}_j \right] = 0
 \end{equation}
 %*

\section{A generalized translation associated with the KdV hierarchy}

In this section we develop the technique of a generalized
translation from \cite{Gordon}.

For $\eta \in \Real$ define an operator $D^\eta$ acting on the
space of functions of one variable as $(D^{\eta} f)(\xi) =
\frac{\xi \eta}{2(\xi-\eta)}f(\eta)$. Define the operator
$\mathcal{B}$ by the rule
 %*
 \begin{equation*}
\mathcal{B}(f,h)(\xi,\eta) = \frac{\xi\,\eta}{2(\xi - \eta)}
(f(\xi) h(\eta) - f(\eta) h(\xi)) = f(\xi) (D^\eta h)(\xi) -
g(\xi) (D^\eta f)(\xi).
 \end{equation*}
 %*
It possesses the following properties:
 \begin{itemize}
\item $ \mathcal{B} ( f, h ) (\xi,\eta)= -\mathcal{B} (h,f)
(\xi,\eta)$.

\item $ \mathcal{B} ( f, h ) (\xi,\eta) = \mathcal{B} ( f, h )
(\eta,\xi) $.

\item $ \mathcal{B} $ is a bilinear operator.

\item If $f(\xi)$ and $h(\xi)$  are polynomials, then $
\mathcal{B} ( f, h ) (\xi,\eta)$ is also a polynomial.

\item $ \mathcal{B} ( f, \xi^{-1} ) (\xi,\eta) = \frac{
\displaystyle f(\xi) \xi - f(\eta) \eta}{ \displaystyle 2
(\xi-\eta) } $.

\item $ \mathcal{B}(1,2 \xi^{-1}) = 1$.
 \end{itemize}

Define also an operator $\mathcal{B}_k$ acting on the set of
$k$-tuples of functions of one variable as follows:
 %*
 \begin{multline*}
\mathcal{B}_k(f_1,\ldots,f_k)(\xi_1,\ldots,\xi_k) =
\frac{\prod_{i=1}^k \xi_i^{k-1}}{2^{k-1} \prod_{1\le i<j\le
k}(\xi_i - \xi_j)}\left(  \sum_{\sigma \in S_k } (-1)^\sigma
f_1(\xi_{\sigma(1)})\ldots f_k(\xi_{\sigma(k)}) \right) =\\
= \frac{\prod_{i=1}^k \xi_i^{k-1}}{
 2^{k-1} W(\xi_1,\xi_2,\ldots,\xi_k)}\, \left|
\begin{array}{cccc} f_1(\xi_1) & f_2(\xi_1) & \ldots &
f_k (\xi_1) \\
f_1(\xi_2) & f_2(\xi_2) & \ldots &
f_k (\xi_2) \\
\vdots & \vdots &  & \vdots \\
f_1(\xi_k) & f_2(\xi_k) & \ldots & f_k (\xi_k)
\end{array}
  \right|
 \end{multline*}
 %*
where $W(\xi_1,\ldots,\xi_k)$ is the Vandermonde determinant.

Note that $\mathcal{B}_1(f) = f$, $\mathcal{B}_2(f,g) =
\mathcal{B} (f,g)$.

Let $f_i$ be a function  of variables $\xi; \,t_1,\dots, t_g$. So
that one has
 %*
 \begin{equation*}
\partial_j \mathcal{B}_k(f_1,\ldots,f_k)(\xi_1,\ldots,\xi_k) =
\sum_{i=1}^k \mathcal{B}_k (f_1,\dots, \partial_j f_i, \dots,f_k)
(\xi_1,\ldots,\xi_k).
 \end{equation*}
 %*

For a fixed function $h$ put, by definition, $(T_\xi^\eta f)(\xi)=
(T(h)_{\xi}^{\eta} f) (\xi) = \mathcal{B}(f,h)(\xi,\eta)$.
 \begin{lem}\label{asss}
The operators $T_{\xi}^{\eta}$ satisfy the associativity condition
$T_{\xi}^{\eta} T_{\xi}^{\tau} = T_{\eta}^{\tau} T_{\xi}^{\eta}$
and the commutativity condition $T_{\xi}^{\eta}
T_{\xi}^{\tau}=T_{\xi}^{\tau} T_{\xi}^{\eta}$.
 \end{lem}

\prf Calculate $T_{\xi}^{\tau} T_{\xi}^{\eta} f$:
 %*
 \begin{multline*}
T_\xi^\tau T_\xi^\eta f  = \frac{\xi \tau}{2(\xi - \tau)} \left(
\xi \eta \frac{ f(\xi) h(\eta) - f(\eta) h(\xi)}{2(\xi-\eta)}\,
h(\tau) - \tau \eta \frac{ f(\tau) h (\eta) - f (\eta) h (\tau)}{
2 (\tau - \eta)} h
(\xi) \right) = \\
\xi^2 \eta^2 \tau^2 \frac{f(\xi) h(\eta) h(\tau)
(\tau^{-1}-\eta^{-1}) + f(\eta) h(\tau) h (\xi)
(\xi^{-1}-\tau^{-1}) +  f(\tau) h(\xi) h(\eta)
(\eta^{-1}-\xi^{-1})}{ 4 (\xi-\tau)(\xi-\eta)(\eta-\tau) } = \\
\mathcal{B}_3 ( f(\xi),  h(\xi), h(\xi)\xi^{-1}) ( \xi,\eta,\tau).
 \end{multline*}
 %*
This expression is invariant under all the permutations of the
variables $\xi$, $\eta$, $\tau$. Lemma is proved. \qed

 \begin{cor}
The operator $T_{\xi}^{\eta}$ is a commutative operator of
generalized translation and
 %*
 \begin{equation*}
T_\xi^\eta \,1 = \frac{\xi \eta}{2 (\xi-\eta)} (h(\eta)-h(\xi)).
 \end{equation*}
 %*
In particular, $T_{\xi}^{\eta} \,1 =1$ if and only if $h(\xi) =
2/\xi$.
 \end{cor}

 \begin{rem}
The generalized translation operator $\mathcal{D}_\xi^\eta (f) =
\frac{\xi f(\xi) -\eta f(\eta)}{\xi-\eta}$ from \cite{Gordon} is
equal to $T^\eta_\xi$ when $h=2/\xi$.
 \end{rem}

 \begin{rem}
Let $h(\xi) = h_{-1}/\xi + \widetilde{h}(\xi)$ where
$\widetilde{h}(\xi)$ is a function regular in a neighbourhood of
the origin. Then for a function $f(\xi)$ regular in the vicinity
of the origin the function $f(\xi,\eta) = T_{\xi}^{\eta} f$ is
regular in the vicinity of the point $(\xi,\eta) = (0,0)$.
 \end{rem}

 \begin{defn}\label{defnpol}
A polarization of a smooth function $f(\xi)$ is a symmetric
function of two variables $f(\xi,\eta)$ such that $f(\xi,\xi) = 2
f(\xi)$.
 \end{defn}

 \begin{lem}\label{polard}
Let $f(\xi,\eta)$ be a polarization of a function $f(\xi)$. Then
 %*
 \begin{equation}\label{polareq}
\left.\frac{\partial f(\xi,\eta)}{\partial \xi}\right|_{\xi=\eta}=
\frac{\partial  f(\xi)}{\partial \xi}.
 \end{equation}
 %*
 \end{lem}
\prf For a symmetric function $f(\xi,\eta)$ there exists a
function $h(s_1, s_2)$ such that
$f(\xi,\eta)=h(\xi+\eta,\xi\eta)$. Since $2f(\xi)=f(\xi,\xi)=h(2
\xi, \xi^2)$, one has
 %*
 \begin{equation*}
\left.\frac{\partial f(\xi,\eta) }{\partial \xi}\right|_{\xi=\eta}
= \left.\frac{\partial h}{\partial s_1}(\xi + \eta, \xi\eta) +
\frac{\partial h}{\partial s_2}(\xi + \eta, \xi\eta) \eta
\right|_{\xi=\eta} = \frac{\partial h}{\partial s_1}(2\xi, \xi^2)
+ \frac{\partial h}{\partial s_2}(2\xi,\xi^2) \xi.
 \end{equation*}
 %*
On the other hand,
 %*
 \begin{equation*}
\frac{\partial f}{\partial \xi}(\xi) = \frac{1}{2} \frac{\partial
h(2 \xi, \xi^2) }{\partial \xi} = \frac{ \partial h }{\partial
s_1}(2 \xi, \xi^2) + \frac{\partial h}{\partial s_2}(2 \xi, \xi^2)
\xi = \left.\frac{\partial f(\xi,\eta) }{\partial
\xi}\right|_{\xi=\eta}\qed
 \end{equation*}
 %*

 \begin{ex}
Let $f(\xi)=\sum_i g_i(\xi) h_i(\xi)$. Then the function
$f(\xi,\eta) = \sum_i \left(g_i(\xi) h_i(\eta) +  g_i(\eta)
h_i(\xi)\right)$ is the polarization of $f(\xi)$.
 \end{ex}

Let $F_n$ is a set of smooth functions on $n$ variables.

 \begin{defn}
Let $G: F_1^k \to F_1$ and $\widehat{G}: F_1^k \to F_2$. The
operator $\widehat{G}$ is called a polarization of the operator
$G$ if the function $\widehat{G}(f_1,\ldots ,f_k)$ is a
polarization of the function $G(f_1,\ldots ,f_k)$ for any
$f_1,\ldots ,f_k$.
 \end{defn}

Recall (\cite{hirota}) that the one-variable Hirota operator
$\mathrm{H}_\xi$ is given by
 %*
 \begin{equation*}
\mathrm{H}_\xi[f(\xi),g(\xi)] = f'(\xi) g(\xi) - f(\xi) g'(\xi).
 \end{equation*}
 %*

 \begin{lem}

The operator $\frac{2}{\xi \eta}\mathcal{B}(f,g)(\xi,\eta)$ gives
the polarization of the Hirota operator
\end{lem}

\prf We need to prove that
 %*
 \begin{equation*}
\lim_{\eta \to \xi} \mathcal{B}(f,g)(\xi,\eta) =
\frac{\xi^2}{2}\mathrm{H}_{\xi}[f(\xi),g(\xi)].
 \end{equation*}
 %*
  Let $\eta=\xi+\varepsilon$. Then $f(\eta) = f(\xi) +
\varepsilon f'(\xi) + O(\varepsilon^2)$ and $g(\eta) = g(\xi) +
\varepsilon g'(\xi) + O(\varepsilon^2)$. Hence
$\mathcal{B}(f,g)(\xi,\xi+\varepsilon) = \frac{\xi
(\xi+\varepsilon)}{-2 \varepsilon} \left( f(\xi) g'(\xi)
\varepsilon - g(\xi) f'(\xi) \varepsilon + O(\varepsilon^2)\right)
= \frac{\xi^2}{2}\mathrm{H}_{\xi}[f(\xi),g(\xi)] +
O(\varepsilon)$\qed

Define operators $D_i$  with the help of the expansion
 %*
 \begin{equation*}
(D^{\eta} f) (\xi) = \sum_{i \in \mathbb{Z}} (D_i f)(\xi) \eta^i.
 \end{equation*}
 %*

\begin{lem}
\label{lemD}
 Let $f(\xi) = \dots + f_0 + f_1 \xi + f_2 \xi^2 + \ldots$.
Then
 %*
 \begin{equation*}
(D_k f)(\xi) = \frac12 (\dots + f_0 \xi^{-k+1} + \dots + f_{k-1}).
 \end{equation*}
 %*

If $f(\xi) = f_0 + f_1 \xi + \dots+f_n \xi^n$ is a polynomial then
$(D_1 f)(\xi) = \frac12 f_0$ and $(D_{n+1} f)(\xi) = \frac12
\,\xi^{-n} f(\xi)$.
 \end{lem}
It is clear that
\begin{equation}
\label{Dk+1f}
 (D_{k+1} f)(\xi) =  \xi^{-1} (D_k f)(\xi) + \frac12
f_k.
 \end{equation}

Note one more property of the operators $D_i$:
 \begin{lem} Let $f(\xi)$ be a polynomial. Then
 %*
 \begin{multline*}
\sum_{k\ge 0, m\ge 0} D_{k+m+1}(f(\xi)) \eta^k \zeta^m  =
\frac{\xi}{2(\eta-\zeta)}\left( \frac{\eta}{\xi-\eta}f(\eta)  -
\frac{ \zeta}{\xi-\zeta}f(\zeta) \right)= \frac{2}{\eta\,\zeta}
\mathcal{B}\left(1, D_{\xi}^\eta f(\xi)\right)(\eta,\zeta).
 \end{multline*}
 %*
 \end{lem}

 Define the operators $d_i$ by the formula $T_\xi^\eta f(\xi)
= \sum (d_i f(\xi)) \eta^i$. Then
 %*
 \begin{equation*}
d_i f(\xi) = f(\xi) D_i h(\xi) - h(\xi) D_i f(\xi).
 \end{equation*}
 %*

From Lemma \ref{asss} using the standard methods we obtain:

 \begin{lem}
The linear space spanned by the operators $d_i$, $i = 1,\dots$, is
an associative and commutative algebra with the following
multiplication:
 %*
 \begin{equation*}
d_i d_j = \sum_{k=0}^{i+j} c_{ij}^k d_k
 \end{equation*}
 %*
where the structure constants $c_{ij}^k $ are found from the
expansion
 %*
 \begin{equation*}
T_\xi^\eta \xi^k = \frac{\xi \eta}{2 (\xi -\eta)} \left( \xi^k
h(\eta) - \eta^k h(\xi) \right) = \sum _{i+j \ge k} c_{ij}^k \xi^i
\eta^j.
 \end{equation*}
 %*
\end{lem}

%\section{The relations between  generating functions}
For the sequence of function $\{u_1,\ldots,u_g\}$ of variables
$t_1 = x, t_2, \dots, t_g$ introduce the generating functions
 %*
 \begin{equation*}
\mathbf{u}(\xi)=\sum_{i=1}^{g} u_{i} \xi^i,\quad
\mathbf{u}'(\xi)=\sum_{i=1}^{g} u_{i}' \xi^i,\quad \dots,
\mathbf{u}^{(k)}(\xi)=\sum_{i=1}^{g} u_{i}^{(k)} \xi^i
 \end{equation*}
 %*
(prime here, as usual, means a differentiation with respect to
$x$). The following statement gives an expression of the third
derivatives $u_1''',\ldots, u_g'''$ in terms of functions
$u_1,\ldots,u_g$ and their first derivatives. Moreover, it allows
to express  these derivatives recursively as a differential
polynomial in $u_1$. This is one of the key results of the paper:

 \begin{thm}
The sequence $\{u_1,u_2,\ldots,u_{g}\}$ satisfies condition
\eqref{4qi+1} of Lemma \ref{lem2} if and only if the generating
function $\mathbf{u}(\xi)$ is a solution of the following
equation:
 %*
 \begin{equation}\label{4mu'}
\mathbf{u}'''(\xi)+2 u_1'(2-\mathbf{u}(\xi)) -
4(\xi^{-1}+u_1)\,\mathbf{u}'(\xi) =0.
 \end{equation}
 %*
 \end{thm}
\prf We have
 %*
 \begin{equation*}
\mathbf{u}'''(\xi)+2 u_1' (2-\mathbf{u}(\xi)) -
4(\xi^{-1}+u_1)\,\mathbf{u}'(\xi)= \sum_{i=1}^g (u_i'''-2
u_1'u_i-4u_1 u_i'-4 u_{i+1}') \xi^i.
 \end{equation*}
 %*
The coefficients at $\xi^i$ in the right-hand side of this formula
are all zero if and only if condition \eqref{4qi+1} of Lemma
\ref{lem2} holds. \qed

%\section{-----------------------------------}
Take, by definition,
 %*
 \begin{equation*}
\partial_k \mathbf{u}(\xi)=\sum_{i=1}^{g} \partial_k u_{i} \xi^i.
 \end{equation*}
 %*

 \begin{lem}
Equations  \eqref{diqj-djqi} and \eqref{di+1} together are
equivalent to the following equation:
 %*
 \begin{equation*}
\partial_{k+1} \mathbf{u}(\xi)= \xi^{-1} \partial_k \mathbf{u}(\xi) -
\frac{1}{2} u_k \mathbf{u}'(\xi) + \frac{1}{2} u'_k
(\mathbf{u}(\xi))-u_k'.
 \end{equation*}
 %*
This equation allows to determine recursively the partial
derivatives $\partial_k \mathbf{u}(\xi)$:
 %*
 \begin{equation}\label{dkq}
\partial_k \mathbf{u}(\xi) = D_k(2-\mathbf{u}(\xi))\mathbf{u}'(\xi) -
D_k(\mathbf{u}'(\xi))(2-\mathbf{u}(\xi)).
 \end{equation}
 %*
 \end{lem}

In the sequel we suppose that \eqref{4mu'} and \eqref{dkq} hold
for the function $\mathbf{u}(\xi)$.

\begin{cor}
 %*
 \begin{eqnarray}
&&\partial_k \mathbf{u}'(\xi) =
D_k(2-\mathbf{u}(\xi))\mathbf{u}''(\xi) -
D_k(\mathbf{u}''(\xi))(2-\mathbf{u}(\xi)),\label{dku'}\\
&&\begin{aligned}
 \partial_k \mathbf{u}''(\xi) &= 4 (\xi^{-1}+u ) \partial_k \mathbf{u}(\xi)
 -2 u_k' (2-\mathbf{u}(\xi)) \\
 & + D_k(\mathbf{u}''(\xi)) \mathbf{u}'(\xi) -
 D_k(\mathbf{u}'(\xi))\mathbf{u}''(\xi).\label{dku''}
\end{aligned}
 \end{eqnarray}
 %*
\end{cor}
Let $\partial(\eta) = \sum_{i=1}^{g} \eta^{i} \partial_i$. Note
that for a fixed $\eta$ the operator $\partial(\eta)$ is an
operator of differentiation in the direction of the vector
$(\eta,\eta^2,\dots,\eta^g)$, i.e.\
 %*
 \begin{equation*}
\partial(\eta) f(t_1,\ldots,t_g) = \left.\frac{ \partial }
{\partial \tau} f(t_1 + \tau\eta, \ldots, t_g +
\tau\eta^g)\right|_{\tau = 0}.
 \end{equation*}
 %*

\begin{cor}
 %*
 \begin{equation*}
\bigl( \partial_x^2 \partial(\xi) + 2(2-\mathbf{u}(\xi))\partial_x
- 4(\xi^{-1}+u_1)\partial(\xi)\bigr) u = 0.
 \end{equation*}
 %*
\end{cor}

\prf Recall that $\partial_i u_1 = u_i'$, and therefore
$\mathbf{u}'(\xi) =
\partial(\xi) u_1$ and $\mathbf{u}'''(\xi) = \partial_x^2\partial(\xi) u$.
The statement follows now from \eqref{4mu'}. \qed

Denote $\mathcal{T}_\xi^\eta = T(2-\mathbf{u}(\xi))_\xi^\eta$. The
operator $\mathcal{T}_\xi^\eta$ plays a special role below, as
shows

 \begin{thm}\label{deuxthm}
 %*
 \begin{eqnarray}
&&\partial(\eta)\mathbf{u}(\xi) =\mathcal{T}_\xi^\eta \partial_x
\mathbf{u}(\xi). \label{deux}\\
&&\partial(\eta) \mathbf{u}'(\xi) = \mathcal{T}_\xi^\eta
\partial_x^2
\mathbf{u}(\xi).\label{deux'}\\
&&\partial(\eta) \mathbf{u}''(\xi) =\mathcal{T}_\xi^\eta
\partial_x^3 \mathbf{u}(\xi)- \mathcal{B}_\xi^\eta (
\mathbf{u}'(\xi), \mathbf{u}''(\xi)).\label{deux''}
 \end{eqnarray}
 %*
\end{thm}

\prf These formulas follow from the definition of the operator
$\mathcal{T}_\xi^\eta$ and equations \eqref{dkq}, \eqref{dku'},
\eqref{dku''}.\qed

Note that \eqref{deux}, \eqref{deux'} imply that
 %*
 \begin{equation*}
[\partial_x, \mathcal{T}_\xi^\eta] \partial_x \mathbf{u}(\xi) =0.
 \end{equation*}
 %*

The associativity condition for the operator
$\mathcal{T}_\xi^\eta$ is equivalent to the following relation,
which will be used later:
 %*
 \begin{equation}\label{dzxideta}
\partial(\zeta) \frac{\xi\,\eta}
{2(\xi-\eta)}\frac{2-\mathbf{u}(\eta)}{2-\mathbf{u}(\xi)} =
\frac{1}{(2-\mathbf{u}(\xi))^2} \mathcal{T}_\xi^\eta
\mathcal{T}_\xi^\zeta \mathbf{u}'(\xi).
 \end{equation}
 %*

%Also we have a nice formula for generating function of the
%operators $\mathcal{U}_i$ :
Now we  describe the family of differential operators
$\{\mathcal{U}_i\}$ using  the method of generating function.

 \begin{lem}
The generating function of the sequence of operators
$\mathcal{U}_i$ is:
 %*
 \begin{equation*}
\sum_{i=1}^g \mathcal{U}_i \xi^i = \frac{1}{2}\bigl( (\mathcal{L}
-\xi^{-1})\partial(\xi) +
\partial(\xi)(\mathcal{L}-\xi^{-1})\bigr) +\frac{1}{4}
\bigl((2-\mathbf{u}(\xi))\partial_x +\partial_x
(2-\mathbf{u}(\xi))\bigr).
 \end{equation*}
 %*
 \end{lem}

\section{The hyperelliptic curve associated with a solution of
KdV}

\begin{thm}
Suppose the generating function $\mathbf{u}(\xi)$ satisfies
\eqref{4mu'} and \eqref{dkq}. Let
 %*
\begin{equation}\label{4mu}
4 \mu(\xi) =  \mathbf{u}'(\xi)^2 + 2 \mathbf{u}''(\xi)\left(2-
\mathbf{u}(\xi)\right)+ 4(\xi^{-1}+u_1)\left(2-
\mathbf{u}(\xi)\right)^2 .
 \end{equation}
 %*
Then $\mu(\xi)=4\xi^{-1}+\sum_{i=1}^{2g}\mu_i \xi^i$ where $\mu_i$
are constants, $i=1,\ldots,2g$.
 \end{thm}

\prf It follows from  \eqref{deux}, \eqref{deux'}, \eqref{deux''}
that $\partial(\eta) \mu(\xi) = 2 \mathbf{u}'(\xi)
\partial(\eta)\mathbf{u}'(\xi) + 2(2-\mathbf{u}(\xi)) \partial(\eta)
\mathbf{u}''(\xi) - 2\mathbf{u}''(\xi) \partial(\eta)
\mathbf{u}(\xi) + 4\mathbf{u}'(\eta) (2-\mathbf{u}(\xi))^2 -
8(\xi^{-1}+u_1) (2-\mathbf{u}(\xi)) \partial(\eta) \mathbf{u}(\xi)
= 0$. Therefore $\partial_i \mu_j =0 $ where $1\le i \le g$, \,
$1\le j \le 2g$, and all the $\mu_i$ are constants.\qed

Assume $u_k=0$ for $k > g$.  Equation \eqref{4mu} implies that
 %*
\begin{equation}\label{uk+1}
u_{k+1} = \frac{1}{4} \mu_k + J_k
(u,u',u'',\ldots,u_{k},u'_{k},u''_{k}),\quad k=1,\dots,2g,
 \end{equation}
 %*
where $J_k$ are polynomials. We see that the functions $u_k$,
$k=2,\ldots,g$, can be expressed recursively via the function
$u_1$, its derivatives and the constants $\mu_i$, namely
 %*
 \begin{equation}\label{4.3}
u_k = \Theta_k(u,u',\ldots, u^{(2k-2)},\mu_1,\ldots,\mu_{k-1})
 \end{equation}
 %*
where $\Theta_k$ are polynomials.

Note that the condition $J'_g = 0$ is equivalent to the stationary
$g$-KdV equation. For $k>g$ one has $J_k=-1/4 \mu_k$, which gives
rise to integrals of the higher KdV equation (see \cite{ZMNP}).

Since $\partial_k u_1=u_k'$, the partial derivative of $u_1$ with
respect to $t_k$ can also be expressed in terms of derivatives
with respect to $x$. Therefore the behavior of function $u_1$
along the coordinate axes $t_2,\ldots,t_g$ can be reconstructed if
its derivatives with respect to $x$ are known.

 \begin{lem}
Let $\mu_i$ be constants. Then equation \eqref{4mu} implies
equation \eqref{4mu'}. If equations \eqref{diqj-djqi} and
\eqref{4mu} hold, then equation \eqref{dkq} holds, too.
 \end{lem}

\prf The first statement of the lemma is trivial. From the
equality $ u_{k+1}'=1/4(u_k'''-2 u' u_k - 4 u u_k')$ one obtains
 %*
 \begin{multline*}
\partial_{k+1} u_m' = \partial_m u_{k+1}' = 1/4 ( \partial_m u_k ''' - 2
u_m'' u_k - 2 u' \partial_m u_k - 4 u_m' u_k' - 4 u \partial_m
u_k') = \partial_k u_{m+1}'+1/2 u_k'' u_m - 1/2 u_m'' u_k.
 \end{multline*}
 %*
This equation proves that \eqref{dku'} holds. So \eqref{dku''} is
also true. Integration of \eqref{dku'} with respect to  $x$ gives
the formula $\partial_k \mathbf{u}(\xi) =
D_k(2-\mathbf{u}(\xi))\mathbf{u}'(\xi) -
D_k(\mathbf{u}'(\xi))(2-\mathbf{u}(\xi))+\varphi(t_2,\ldots,t_g)$
where the function $\varphi$ does not depend on $x$. Combining the
last equation with \eqref{4mu}, we obtain
 %*
 \begin{multline*}
0=4\, \partial_k \mu(\xi)= \partial_k (\mathbf{u}'(\xi)^2 + 2
\mathbf{u}''(\xi)\left(2- \mathbf{u}(\xi)\right)+
4(\xi^{-1}+u_1)\left(2-
\mathbf{u}(\xi)\right)^2) =\\
\left(8 (\xi^{-1}+u_1)(\mathbf{u}(\xi)) + 2 \mathbf{u}''(\xi)
\right)\varphi(t_2,\ldots,t_{g}).
 \end{multline*}
 %*
The function in parenthesis cannot vanish identically as a
function of $x$, and thus $\varphi(t_2,\ldots,t_{g}) \equiv 0$.
\qed

Summarize the results obtained:
\begin{thm} \label{main1}
The following statements are equivalent:
\begin{enumerate}
\item The function  $u_1$ is a solution of a stationary $g$-KdV
equation. \item There exists a sequence of functions
$\{u_1,\ldots,u_g\}$ such that the operators
$\mathcal{L}=\partial_x^2 -u_1$ and $\mathcal{U}_i=\partial_x ^2
\partial_i -\frac{1}{2}( u_1 \partial_i + \partial_i u_1)-
\frac{1}{4} ( u_i \partial_x +  \partial_x u_i ) -\partial_{i+1}$,
$1 \le i \le g$ commute. \item There exists a sequence of
functions $\{u_1,\ldots,u_g\}$ and a set of constants
$\mu_1,\ldots,\mu_{2g}$ such that the  generating function
$\mathbf{u}(\xi)=\sum_{i=1}^g u_{i} \xi^i$ satisfies \eqref{4mu}.
\end{enumerate}
\end{thm}

In order to find out the relation between the constants $\mu_k$
and the coefficients $a_i$ we need the following result:

 \begin{lem}\label{lem15}
The operator $\mathbf{U}=\mathcal{U}_1
\mathcal{L}^{g-1}+\mathcal{U}_2 \mathcal{L}^{g-2}+\cdots
+\mathcal{U}_{g}$
\begin{enumerate}
\item commutes with the operator $\mathcal{L}$;

\item is an operator of order $2g+2$ with the leading coefficient
$1$;

\item contains the differentiation with respect to $x$ only.
\end{enumerate}
 \end{lem}

\prf The first statement of the lemma is obvious. The leading term
of $\mathbf{U}$ is a composition of the leading terms of the
operators $\mathcal{U}_1$ and $\mathcal{L}^{g-1}$, so it is equal
to $\partial_x^{2g+2}$. This proves the second statement. Since
$\mathcal{U}_i = \partial_i \mathcal{L} -\partial_{i+1} - 1/2 u_i
\partial_x +1/4 u_i'$, the sum $\mathcal{U}_{i} \mathcal{L} -
\mathcal{U}_{i+1}$ does not contain the differentiation with
respect to $t_{i+1}$. By recursion we get the third statement of
the lemma. \qed

\begin{thm}\label{thmUL}
Under the hypotheses of Theorem \ref{main1} the operator $A$ can
be decomposed as $A =\mathcal{U}_1 \mathcal{L}^{g-1}+\mathcal{U}_2
\mathcal{L}^{g-2}+\ldots +\mathcal{U}_{g}$, where $A = A_g +
\sum_{i=0}^{g-2} a_i A_i$. The coefficients $\mu_k$ and $a_i$
satisfy the following relation:
 %*
 \begin{equation}\label{mukaij}
\mu_{k}=8 a_{g-k-1}+ 4 \sum_{i=1}^{k-2}a_{g-i-1} a_{g-k+i}.
 \end{equation}
 %*
for $k=1,\ldots,g-1$.
\end{thm}

\prf The first statement of the theorem follows from lemma
\ref{lem15} and uniqueness of the operators $A_k$ (see Theorem
\ref{thmAbasis}).

The function $u_k$ is a differential polynomial $u_k =
\Theta_k(u_1,\;u_1',\ldots,u_1^{(2k)},\mu_1,\ldots,\mu_{k-1})$.
Let $\varepsilon_k$ be a constant term of $\Theta_k$. Then
$\mathcal{U}_1 \mathcal{L}^{g-1}+\mathcal{U}_2
\mathcal{L}^{g-2}+\cdots +\mathcal{U}_g = \partial_x^{2g+1} - 1/2
\sum u_k \partial_x^{2g-2k-1} + \sum \vartheta_i \partial_x^i =
\partial_x^{2g+1} - 1/2 \sum \varepsilon_k \partial_x^{2g-2k-1} + \sum
\widetilde{\vartheta_i} \partial_x^i$ where $\vartheta_i$ and
$\widetilde{\vartheta_i}$ are differential polynomials in $u_1$
without constant terms. On the other hand, $A = \partial_x^{2g+1}
- 1/2 \sum a_k
\partial_x^{2k+1} + \sum \widetilde{\vartheta_i} \partial_x^i$. Thus, $ a_k
= -1/2 \varepsilon_{g-k-1}$, and so it remains to find
$\varepsilon_{k}$. The result follows now from \eqref{4mu}.\qed

%In such a way we obtain one of the main results.
The following corollary is one of the main results of the paper:

\begin{cor}\label{cor12}

There is a canonical way to associate a solution $u_1$ of the
stationary $g$-KdV equation with a hyperelliptic curve
 %*
\begin{equation}\label{Gammadef}
\Gamma=\{(\xi,y) \in \mathbb{C}^2\;|\; y^2=4 \mu(\xi)\}.
 \end{equation}
 %*
The coefficients $\mu_1,\ldots,\mu_{g-1}$ are expressed in terms
of the constants $a_i$ as in equation \ref{mukaij}, and
$\mu_g,\ldots,\mu_{2g}$ are found from \ref{4.3} in terms of  the
values of $u_1^{(k)}(t_0)$, $k=0,1,\ldots$,  at some point $t_0
\in \mathbb{C}^g$.
\end{cor}

 \begin{rem}
The hyperelliptic curve constructed above coincides with the
spectral curve introduced in \cite{3} when the solution $u_1$ is
periodic as a function of $x$. Our construction uses only the
local properties of the function $u_1$, while in \cite{3} only
periodical or rapidly decreasing functions are discussed.
 \end{rem}

 \begin{rem}
The number of singular points on $\Gamma$ is an important
characteristic of the solution $u_1$. This number can be expressed
in terms of $u_1^{(k)}(t_0)$ using the resultant.
 \end{rem}

\section{Fiber bundles associated with the stationary $g$-KdV equations}
\label{sectionR}

The equations described by \eqref{novikov}  are ordinary
differential equations of order $2g+1$, and so their solution
$u_1$ is \textit{uniquely} determinated in a neighbourhood of a
given point $x_0$ by the values $c_k = u_1^{(k)}(x_0)$,
$k=0,\ldots,2g$. Since the coefficients of the KdV equations are
constants, we can take $x_0=0$. The stationary $g$-KdV equations
depend on the numbers $a_{0},\ldots,a_{g-2}$ and so the space of
all such equations is isomorphic to $\mathbb{C}^{g-1}$.

The space  $\mathcal{M}_g$ of all the hyperelliptic curves
$\Gamma=\{(\xi,y) \in \mathbb{C}^2 \mid y^2=4 \mu(\xi)\}$ can be
parametrized by the numbers $\mu_1,\ldots\mu_{2g}$, so it is
isomorphic to $\mathbb{C}^{2g}$.

Denote by $\mathbf{R}_g$ the space of solutions $u$ of all the
stationary $g$-KdV equations such that $u$ is regular at the point
$x_0$. As it was explained above, we can identify the space
$\mathbf{R}_g$ with $\mathbb{C}^{3g}$ using coordinates $(c_0,
c_1, \ldots, c_{2g}, a_0, \ldots, a_{g-2})$.

There exists a canonical map $\pi_{\mathcal{M}}: \mathbf{R}_g \to
\mathcal{M}_g$, which sends a solution $u$ to a hyperelliptic
curve $\Gamma$ described in \eqref{Gammadef}.

Denote by $\mathfrak{U}_g$ the space of $g$-th symmetric powers of
hyperelliptic genus $g$ curves. We consider the  universal bundle
$(\mathfrak{U}_g,\, \mathcal{M}_g,\,  \pi_{\mathfrak{U}})$ where
the  natural projection  $ \pi_{\mathfrak{U}}: \mathfrak{U}_g \to
\mathcal{M}_g $ is given by $ \pi_{\mathfrak{U}_g}( x \in Sym^g
\Gamma ) = \Gamma $.

 \begin{thm}
There exists a canonical fiber-preserving birational equivalence
$\mathbf{R}_g \to \mathcal{U}_g$.
 \end{thm}
\prf Let $\Gamma$ be a hyperelliptic curve associated with the
solution $u_1$ of the stationary $g$-KdV equation (see Corollary
\ref{cor12}). Let $\xi_1,\ldots,\xi_g$ be the roots of the
equation $2-\mathbf{u}(0,\xi)=0$. Denote
$y_i=\mathbf{u}'(0,\xi_i)$. Equation \eqref{4mu} implies that
$y_i^2 = \mathbf{u}'(0,\xi_i)^2 = 4\mu(\xi)$, so the point
$(\xi_i,y_i)$ belongs to $\Gamma$. Thus we have a map $\upsilon:
\mathbf{R}_g \to \mathfrak{U}_g$ given by the formula $\upsilon
(u_1) = (\Gamma,[(\xi_1,y_1),\ldots,(\xi_g,y_g)])$ where
$(\xi_i,y_i)\in \Gamma$. Apparently, $\upsilon$ is
fiber-preserving.

On the other hand, if a curve $\Gamma$ and a point
$[(\xi_1,y_1),\ldots,(\xi_g,y_g)] \in Sym^g \Gamma$ are given,
then in the case of distinct points $(\xi_1,\ldots,\xi_g)$, it is
possible to construct the point
$(c_0,\ldots,c_{2g},a_0,\ldots,a_{g-2})$ as follows. The constants
$a_k$ are a solution of \eqref{mukaij} where the parameters
$\mu_i$ are known. The values $u_i(0)$ are  the symmetric
functions of $\xi_1,\ldots,\xi_g$, namely $u_i(0) = 2 \,
(-1)^{g-i}\, \sigma_{g-i+1}(\xi_1,\ldots,\xi_g)/\sigma_g
(\xi_1,\ldots,\xi_g)$. Then the values $u_i'(0)$ can be found, as
the coefficients of the generating function $\mathbf{u}'(0,\xi)$,
from the equations $\mathbf{u}'(0,\xi_i)=y_i$. All the higher
derivatives $c_k=u_1^{(k)}(0)$ can be found by recursion using
equation \eqref{4.3}. Thus the inverse rational map
$\upsilon^{-1}$ is constructed. \qed

In case of the universal bundle of Jacobians over the moduli space
of genus $g$ hyperelliptic curves this theorem gives the famous
results of Dubrovin and Novikov, see \cite{3}.

\section{Algebraic relations between the operators $\mathcal{L},
\mathcal{U}_1, \ldots, \mathcal{U}_{g}$}

The  Burchnall--Chaundy lemma (\cite{burchaundy}) says that two
commuting differential operators of one variable are always
connected by an algebraic relation. In \cite{krich} the case of
commuting differential operators of $n$ variables was considered.
In the same work they indtroduced a class of $n$-algebraic
families of operators, i.e.\ families of commuting operators
characterized by finite-dimensional algebraic manifolds. The
family $\{\mathcal{L},\mathcal{U}_1,\ldots,\mathcal{U}_g\}$ gives
an example of $n$-algebraic operators from \cite{krich}.

 \begin{lem}
The operators $\mathcal{L},\mathcal{U}_1,\ldots,\mathcal{U}_{g}$
satisfy the following algebraic relation:
 %*
 \begin{equation*}
4(\mathcal{U}_1 \mathcal{L}^{g-1}+\mathcal{U}_2
\mathcal{L}^{g-2}+\ldots+\mathcal{U}_{g-1} \mathcal{L} +
\mathcal{U}_{g} )^2 = ( 4 \mathcal{L}^{2g+1} + \mu_1
\mathcal{L}^{2g-1} + \mu_2 \mathcal{L}^{2g-2} + \ldots+ \mu_{2g}
).
 \end{equation*}
 %*
 \end{lem}

Using the notations $\mathcal{U}(z)= \mathcal{U}_1 z^{g-1} +
\mathcal{U}_2 z^g +\ldots+ \mathcal{U}_{g}$ and
$\widetilde{\mu}(z) = 4 z^{2g+1}+ \mu_1 z^{2g-1}+\ldots+\mu_{2g}$,
one can write down this relation as  $4 \mathcal{U}(\mathcal{L})^2
= \widetilde{\mu}(\mathcal{L})$.

\prf Denote
 %*
 \begin{equation*}
S_i = \frac{1}{2}\, u_i \partial_x -\frac{1}{4}\, u_i'.
 \end{equation*}
 %*
Then
 %*
 \begin{equation}\label{Ui=Si}
\mathcal{U}_i=\partial_i \mathcal{L} - S_i -\partial_{i+1}.
 \end{equation}
 %*
We have
 %*
 \begin{equation*}
[\mathcal{L},S_i] = u_i' \mathcal{L} -u'_{i+1},
 \end{equation*}
 %*
which implies the equation
 %*
 \begin{multline*}
\sum_{i+j=k} \mathcal{U}_i \mathcal{U}_j = \sum_{i+j=k}
(\partial_i \mathcal{L} - S_i-\partial_{i+1})(\partial_j
\mathcal{L} - S_j-\partial_{j+1}) =\\
\sum_{i+j=k} \partial_i \partial_j \mathcal{L}^2 - 2 \partial_i
\partial_{j+1} \mathcal{L} + \partial_{i+1} \partial_{j+1} -
(\partial_i S_j + S_i \partial_j) \mathcal{L} + (\partial_{i+1}
S_j +S_i
\partial_{j+1}) + S_i S_j.
 \end{multline*}
 %*

A direct calculation gives that
 %*
 \begin{eqnarray*}
&&S_i S_j = 1/4 u_i u_j\,
\partial_x^2 + 1/8 ( u_i u_j' - u_j u_i') \partial_x
-1/8 u_i u_j'' +1/16 u_i' u_j',\\
&&\partial_x S_i + S_i \partial_x = u_i \partial_x^2 - 1/4 u_i''.
 \end{eqnarray*}
 %*

Therefore,
 %*
 \begin{multline*}
\sum_{1\le i,j \le g} \mathcal{U}_i \mathcal{U}_j
\mathcal{L}^{2g-i-j+1} =
\partial_x^2 \mathcal{L}^{2g} - \sum_{1\le i \le g} (u_{i} \partial_x^2 -
1/4 u_{i}'') \mathcal{L}^{2g-i} +\\
+\sum_{1\le i,j \le {g}} (1/4 u_i u_j \partial_1^2 -1/2 u_i u_j''
+1/16
u_i' u_j')\mathcal{L}^{2g-i-j} =\\
\mathcal{L}^{2g+1}+u_1 \mathcal{L}^{2g} - \sum_{1\le i \le g}
u_{i} \mathcal{L}^{2g-i+1}-\sum_{1\le i \le g} (u_{i} u_1 - 1/4
u_{i}'')
\mathcal{L}^{2g-i-1} +\\
+\sum_{1\le i,j \le g} (1/4 u_i u_j \mathcal{L}^{2g-i-j+1}
)+\sum_{1\le i,j \le g} (1/4 u_i u_j u_1 -1/8 u_i u_j'' +1/16 u_i'
u_j')\mathcal{L}^{2g-i-j}.
 \end{multline*}
 %*

We see that the coefficient at $\mathcal{L}^{2g-i}$ in this
formula is exactly the coefficient at $\xi^i$ in the expression
$1/16 ( \mathbf{u}'(\xi)^2+2 \mathbf{u}''(\xi)
(2-\mathbf{u}(\xi))+4 (u_1+\xi^{-1})(2-\mathbf{u}(\xi))^2 = 1/4
\mu(\xi)$.\qed

 \begin{cor}
Let $\Psi(t_1,t_ 2,\ldots, t_{g})$ be a common eigenfunction of
the operators $\mathcal{L},\mathcal{U}_1,\ldots,\mathcal{U}_{g}$,
with the eigenvalues $E, \alpha_1,\ldots, \alpha_{g}$. Let
$\xi=E^{-1}$ and $\alpha(\xi) = \sum_{i=1}^{g} \alpha_i \xi^{i}$.
Then
 %*
 \begin{equation}\label{ak}
4 \alpha(\xi)^2 = \mu(\xi).
 \end{equation}
 %*
\end{cor}

\section{A common eigenfunction of the family
$\{\mathcal{U}_i\}$}\label{sectionF}

In this section we construct a common eigenfunction of the family
of commuting differential operators $\{\mathcal{U}_i\}$.

 \begin{lem}\label{4}
 %*
 \begin{equation*}
\frac{\partial}{\partial t_i} \left(
\frac{D_j(2-\mathbf{u}(\xi))}{2-\mathbf{u}(\xi)} \right) =
\frac{\partial}{\partial t_j} \left(
\frac{D_i(2-\mathbf{u}(\xi))}{2-\mathbf{u}(\xi)} \right)
 \end{equation*}
 %*
 \end{lem}

\prf It follows from the definition of the operators
$\partial(\eta)$ and $D_i$ that the expression
$\frac{\partial}{\partial t_i} \left(
\frac{D_j(2-\mathbf{u}(\xi))}{2-\mathbf{u}(\xi)} \right)$ equals
to the coefficient at $\zeta^i\eta^j$ in the expansion
$\partial(\zeta) \frac{\xi\,\eta}{2(\xi-\eta)}
\frac{2-\mathbf{u}(\eta)}{2-\mathbf{u}(\xi)} $ with respect to
$\eta$ and $\zeta$. This function is equal to
$\frac{1}{(2-\mathbf{u}(\xi))^2} \mathcal{T}_\xi^\eta
\mathcal{T}_\xi^\zeta \mathbf{u}'(\xi)  $ (see \eqref{dzxideta}).
 Since the generalized translation
$\mathcal{T}_\xi^\eta$ is commutative, this function is symmetric
with respect to the variables $\zeta$ and $\eta$. Consequently the
coefficients of $\zeta^i\eta^j$ and $\zeta^j\eta^i$ are equal.
\qed

 \begin{cor}
There exists a function  $F(\xi)=F(t_1,\ldots,t_g,\xi)$ such that
$\partial_i F =\frac{D_i(2-\mathbf{u}(\xi))}{2-\mathbf{u}(\xi)}$,
$1\le i \le g$. The function  $F(\xi)$ is uniquely defined up to
an additive constant in a neighborhood of any point
$(\bar{t}_0,\xi_0)=(t^0_1,\ldots,t^0_g, \xi_0)$ such that
 $2-\mathbf{u}(\bar t_0;\xi_0)\ne0$.
 \end{cor}

Consider also the function
$\Phi=\Phi(t_1,\ldots,t_g;E=\xi^{-1},\alpha_1,\ldots,\alpha_g)$
given by
 %*
 \begin{eqnarray}\label{Phi}
\Phi= \sqrt{2-\mathbf{u}(\xi)} \exp  \left(2 \alpha(\xi) F(\xi)
\right) \exp (- 2 \sum_{i=1}^{g} D_i(\alpha(\xi)) t_i ),
 \end{eqnarray}

where $\alpha(\xi) =\sum_{i=1}^{g} \alpha_i \xi^{i}$. The function
$\Phi$ is uniquely defined up to a multiplicative constant in a
neighborhood of any point $(\bar t_0,\xi_0,\bar \alpha)$ such that
$2-\mathbf{u}(\bar t_0;\xi_0)\ne0$.

Find the derivatives of the function  $\Phi$ with respect to
$x=t_1$ and $t_k$, $k \ge 2$:
 %*
 \begin{eqnarray}
&&\Phi' = \frac{ 4 \alpha(\xi) -
\mathbf{u}'(\xi)}{2(2-\mathbf{u}(\xi))}
\Phi,\label{Phi'}\\
&&\begin{aligned} \Phi'' &= \left(\frac{
-\mathbf{u}''(\xi)(\mathbf{u}(\xi) + \mathbf{u}'(\xi)(
4\alpha(\xi) -\mathbf{u}'(\xi))}{ 2 (2-\mathbf{u}(\xi)} + \frac{
(4 \alpha(\xi) -\mathbf{u}'(\xi))^2} {4
(2-\mathbf{u}(\xi))^2}\right)\Phi\\
 &= \frac{16 \alpha(\xi)^2 - 2
\mathbf{u}''(\xi)(2-\mathbf{u}(\xi))-\mathbf{u}'(\xi)^2}{4
(2-\mathbf{u}(\xi))^2} \Phi,
\end{aligned}\label{Phi''}\\
&&\partial_k \Phi = \left(\frac{4 \alpha(\xi)
D_k(2-\mathbf{u}(\xi))-\partial_k \mathbf{u}(\xi)}{2
(2-\mathbf{u}(\xi))}-2 D_k(\alpha(\xi)) \right)\Phi.\label{Phik}
 \end{eqnarray}
 %*

 \begin{lem}\label{lem17}
The function $\Phi$ is an eigenfunction of the operator
$\mathcal{L}$ with the eigenvalue $E \ne 0$ if and only if
$\xi=E^{-1}$ and $\{\xi,4 \alpha(\xi)\} \in \Gamma$, where
$\Gamma$ is a curve defined by equation \eqref{Gammadef}.
 \end{lem}

\prf Equation \eqref{4mu} implies that
 %*
 \begin{equation*}
(\mathcal{L} - E)\Phi = \left(\frac{ 4 \alpha(\xi)^2 - \mu(\xi)}
{2-\mathbf{u}(\xi)}-(\xi^{-1}-E)\right)\Phi.
 \end{equation*}
 %*
The function in parentheses vanishes identically iff $4
\alpha(\xi)^2 - \mu(\xi)=(\xi^{-1}-E)(2-\mathbf{u}(\xi))$.
Differentiating the last formula with respect to $x$,  we obtain
that  $(\xi^{-1}-E) \mathbf{u}'(\xi) = 0$. Hence $\xi^{-1} = E$
and $4 \alpha(\xi)^2 = \mu(\xi)$. \qed

 \begin{thm}
Suppose that $4 \alpha(\xi)^2 = \mu(\xi)$. Then
 \begin{enumerate}
\item The function $\Phi$ is a common eigenfunction of the family
$\mathcal{L}, \mathcal{U}_1,\ldots,\mathcal{U}_g$ with eigenvalues
$E=\xi^{-1}, \alpha_1,\ldots,\alpha_g$.

\item The space of common eigenfunction of operators $\mathcal{L},
\mathcal{U}_1,\ldots,\mathcal{U}_g$ with eigenvalues
$E=\xi^{-1},\alpha_1,\ldots,\alpha_g$ is one-dimensional.
 \end{enumerate}
 \end{thm}

\prf Express the operators $\mathcal{U}_k$ as
 %*
 \begin{equation}\label{U=L+l}
\mathcal{U}_k=\partial_k (\partial_x^2 - (u_1+\xi^{-1})) +
\xi^{-1}
\partial_k - 1/2 u_k \partial_x +1/4 u'_k -
\partial_{k+1}.
 \end{equation}
 %*
Let $\Psi$ be a common eigenfunction of $\mathcal{L}$,
$\mathcal{U}_k$ with the eigenvalues mentioned above. Then
 %*
\begin{equation}\label{hrum}
\left(\xi^{-1} \partial_k - 1/2 u_k \partial_x +1/4 u'_k -
\partial_{k+1} \right) \Psi = \alpha_k  \Psi.
 \end{equation}
 %*
This allows to express all the partial derivatives $\partial_{k}
\Psi$ in terms of $\Psi$ and $\Psi'$, namely
 %*
 \begin{equation*}
\partial_{k} \Psi = D_{k}(2-\mathbf{u}(\xi)) \Psi' + \frac{1}{2}
D_{k}(\mathbf{u}'(\xi)) \Psi - 2 D_k(\alpha(\xi))\Psi,\quad 1 \le
k \le g-1.
 \end{equation*}
 %*
For $k=g-1$ one gets from \eqref{hrum}
 %*
 \begin{equation*}
\xi^{-1} \partial_{g-1} \Psi - 1/2 u_{g-1} \Psi' +1/4 u'_{g-1}
\Psi= \alpha_g \Psi.
 \end{equation*}
 %*
Therefore,
 %*
 \begin{equation*}
\xi^{-1} \left (D_{g}(2-\mathbf{u}(\xi)) \Psi' + \frac{1}{2}
D_{g}(\mathbf{u}'(\xi)) \Psi - 2 D_g(\alpha(\xi))\right)- 1/2 u_g
\Psi' +1/4 u'_g \Psi- \alpha_g \Psi = 0.
 \end{equation*}
 %*
Using Lemma \ref{lemD} and \eqref{Dk+1f} we obtain
 %*
 \begin{equation*}
(2-\mathbf{u}(\xi))) \Psi' = (1/2 \mathbf{u}'(\xi)+\alpha(\xi))
\Psi.
 \end{equation*}
 %*
Thus,
 %*
 \begin{equation*}
\frac{\Psi'}{\Psi}= \frac{-1/2
\mathbf{u}'(\xi)+\alpha(\xi)}{2-\mathbf{u}(\xi)},
 \end{equation*}
 %*
and
 %*
 \begin{multline*}
\frac{\partial_k \Psi}{\Psi} = D_k(2-\mathbf{u}(\xi)) \frac{-1/2
\mathbf{u}'(\xi)+\alpha(\xi)} {2-\mathbf{u}(\xi)} + \frac{1}{2}
D_{k}(\mathbf{u}'(\xi)) -2 D_k(\alpha(\xi))\\
= \frac{ \alpha(\xi)D_k(2-\mathbf{u}(\xi)) - 1/2 \partial_k
\mathbf{u}(\xi)}{2-\mathbf{u}(\xi)}- 2 D_k(\alpha(\xi)),\quad
k=2,\ldots,g.
 \end{multline*}
 %*

We see that $\frac{\partial_k \Psi}{\Psi}=\frac{\partial_k
\Phi}{\Phi}$, $k=1,\ldots,g$. Therefore, $\Psi = \lambda  \Phi$
where $\lambda$ is a constant. \qed

Consider now the special case $E=0$.
 \begin{thm}
The space of common eigenfunction of operators $ \mathcal{L},
\mathcal{U}_1,\ldots,\mathcal{U}_g$ with eigenvalues
$0,\alpha_1,\ldots,\alpha_g$, where  $4 \alpha_g^2 = \mu_{2g}$, is
one-dimensional.
 \end{thm}

\prf Let $\Phi^0$ be a common eigenfunction of operators
$\mathcal{L}, \mathcal{U}_1,\ldots,\mathcal{U}_g$ with eigenvalues
$0,\alpha_1,\ldots,\alpha_g$. Since $\mathcal{L}\Phi^0=0$,
equation \eqref{Ui=Si} implies that
 %*
 \begin{equation*}
\mathcal{U}_k \Phi^0 = (-1/2 u_k \partial_x +1/4 u_k' -
\partial_{k+1}) \Phi^0 = \alpha_k \Phi^0.
 \end{equation*}
 %*
Therefore,
 %*
 \begin{equation}\label{dkPhi0}
\partial_k \Phi^0 = -(1/2 u_{k-1} \partial_x -1/4 u_{k-1}'+\alpha_{k-1}
)\Phi^0 ,\quad k=2,\ldots,g
 \end{equation}
 %*
and
 %*
 \begin{equation}\label{dxPhi0}
\partial_x \Phi^0 = \frac{ u_{g}'-4 \alpha_{g}}{2 u_{g}} \Phi^0.
 \end{equation}
 %*
Note that \eqref{dxPhi0} and \eqref{4mu} imply
 %*
 \begin{equation*}
\partial_x^2 \Phi^0  = \frac{ 2 u''_{g} u_{g} - (u_{g}')^2 + 16
\alpha^2_{g} } { 4 u_{g}^2 } \Phi_0 = \frac{ 2 u''_{g} u_{g} -
(u_{g}')^2 + 4 \mu_{2g} } { 4 u_{g}^2 } \Phi^0 = u  \Phi^0.
 \end{equation*}
 %*

It follows from \eqref{dkPhi0} and \eqref{dxPhi0} that
 %*
 \begin{equation*}
\partial_k \Phi^0 = \frac { u_{k-1}( 4 \alpha_{g} - u'_{g}) -u_{g}( 4
\alpha_{k-1} - u'_{k-1})}{4 u_{g}} \Phi^0.
 \end{equation*}
 %*

Since the logarithmic derivatives are uniquely defined, the space
of eigenfunction with eigenvalues  $E=0, \alpha_1, \ldots,
\alpha_{g} = 1/4 \sqrt{\mu_{2g}}$ is one-dimensional. \qed

Note that the function $\Phi$ can be expressed as
 %*
 \begin{equation}\label{phiexpF}
\Phi =\exp\widetilde{F},
 \end{equation}
 %*
where
 %*
 \begin{equation}\label{tildeF}
\widetilde{F} = \left(2\alpha(\xi)F-2\sum_{i=1}^g D_i(\alpha(\xi))
t_i + \frac{1}{2} \log(2-\mathbf{u}(\xi))\right).
 \end{equation}
 %*
We have
 %*
 \begin{equation}\label{ditildeF}
\partial_i \widetilde{F}= \frac{D_i(2-\mathbf{u}(\xi))(4
\alpha(\xi)-\mathbf{u}'(\xi))-D_i(4
\alpha(\xi)-\mathbf{u}'(\xi))(2-\mathbf{u}(\xi))}{2(2-\mathbf{u}(\xi))}.
 \end{equation}
 %*
Using the notation  $\partial(\eta) = \sum_{i=1}^g \eta^i
\partial_i$, we can rewrite these formulas as
 %*
 \begin{multline*}
\partial(\eta)\widetilde{F} =
\frac{\xi\,\eta}{\xi-\eta}\frac{(2-\mathbf{u}(\eta))(4
\alpha(\xi)-\mathbf{u}'(\xi))-(4
\alpha(\eta)-\mathbf{u}'(\eta))(2-\mathbf{u}(\xi))}{4(2-\mathbf{u}(\xi))}
= \frac{\mathcal{T}_\xi^\eta(4 \alpha(\xi)-\mathbf{u}'(\xi))}
{2(2-\mathbf{u}(\xi))}.
 \end{multline*}
 %*

 \begin{lem}
The function $\chi(\xi)= \,\partial_1 \widetilde{F} = \frac{4
\alpha(\xi)-\mathbf{u}'(\xi)}{2(2-\mathbf{u}(\xi))}$ satisfies the
Riccati equation $ \chi'(\xi)+\chi(\xi)^2= u_1+ \xi^{-1}$.
Moreover,
 %*
 \begin{equation*}
\frac{1}{2-\mathbf{u}(\eta)}\,\partial(\eta)\widetilde{F} =
\frac{\xi\,\eta}{2(\xi-\eta)} ( \chi(\xi) - \chi(\eta)).
 \end{equation*}
 %*
 \end{lem}

Denote by $V$ a hypersurface in $\mathbb{C}^{g+1}=
\{(\xi,\alpha_1,\ldots,\alpha_{g})\}$ defined by equation
\eqref{ak}. Recall that $\Gamma=\{\xi, y \in \mathbb{C}^2 \mid y^2
= 4 \mu(\xi) \}$. In coordinates  $E,\alpha_i$ the hypersurface
$V$ is given by the equation
 %*
 \begin{equation*}
4 \left(\alpha_1 E^{g-1}+\alpha_2 E^{g-2}+\ldots+\alpha_{g}
\right)^2 = 4 E^{2g+1}+\mu_1 E^{2g-1}+\ldots+ \mu_{2g},
 \end{equation*}
 %*
and the curve $\Gamma$ is given by the equation $ \eta^2 = 4(4
E^{2g+1}+\mu_1 E^{2g-1}+\ldots+ \mu_{2g})$ where $\eta = y
\xi^{-g}$.

Define a projection $\pi: V \to \Gamma$ by the formula  $\pi (\xi,
\alpha_1, \ldots, \alpha_{g} ) = (\xi,  2 \alpha(\xi) ) $. In the
sequel we will consider the curve $\Gamma$ as subvariety  of $V$
using a canonical embedding $i: \Gamma \hookrightarrow V $ defined
as $i(\xi, \eta) = (\xi, 0, 0, \ldots, \eta)$. Let $V^* = \pi^{-1}
(\Gamma^*)$ where $\Gamma^*=\{ (\xi,y)\in \Gamma;\; \xi \ne 0 \}$.

Recall that $\mathbb{C}^* = \{ \xi \in \mathbb{C} \mid \xi \ne 0
\}$. The function $\Phi$ of equation \eqref{Phi} is defined in the
space $\mathbb{C}^g \times \mathbb{C}^* \times  \mathbb{C}^g$
parametrized with coordinates $ t_1,\ldots,t_g, \xi,
\alpha_1,\ldots,\alpha_g$. Consider this domain as a graded space
using the following grading: $\deg t_k = 1-2k, \; \deg \xi = -2,\;
\deg \alpha_k = 2k+1$. Take also $\deg \mu_i =2i+2$. Then the
equation $4\alpha(\xi)=\mu(\xi)$ defining the variety $V^*$ is
homogeneous.

 \begin{lem}
Let $\Phi$ be a common  eigenfunction of the operators
$\mathcal{L}, \mathcal{U}_1,\ldots,\mathcal{U}_g$ with eigenvalues
$E=\xi^{-1},\alpha_1,\ldots,\alpha_g$. Let
$\gamma_2,\ldots,\gamma_g \in \mathbb{C}$ be arbitrary constants.
Then the function  $\widetilde{\Phi} = \Phi\, \exp (\gamma_2 t_2
+\ldots+ \gamma_g t_g)$ is also a common eigenfunction of these
operators, its eigenvalues given by
$E=\xi^{-1},\;\widetilde{\alpha}_1=\alpha_1-\gamma_2,\;\widetilde{\alpha}_2
=\alpha_2-\gamma_3+\xi \gamma_2,\ldots,
\widetilde{\alpha}_i=\alpha_i-\gamma_{i+1}+\xi \gamma_i, \ldots,
\widetilde{\alpha}_g=\alpha_g+\xi \gamma_g$.
 \end{lem}

\prf Take $\gamma_1=0$. It is obvious that $\frac{\partial_k
\widetilde{\Phi}}{\widetilde{\Phi}} = \frac{\partial_k
\Phi}{\Phi}+\gamma_k$, $k=1,\ldots,g$ and
$\mathcal{L}\widetilde{\Phi}=E \widetilde{\Phi}$. From
\eqref{U=L+l} one obtains that
 %*
 \begin{equation*}
\mathcal{U}_k \widetilde{\Phi} = (\xi^{-1} \partial_k -\frac{1}{2}
u_k
\partial_x + \frac{1}{4} u_k'-\partial_{k+1})\widetilde{\Phi}=
\exp(\gamma_2 t_2+\ldots+\gamma_g t_g) \left(\mathcal{U}_k
\Phi+(\xi^{-1} \gamma_k-\gamma_{k+1})\Phi \right) .
 \end{equation*}
 %*
Therefore,
 %*
\begin{equation*}
\frac{\mathcal{U}_k
\widetilde{\Phi}}{\widetilde{\Phi}}=\widetilde{\alpha}_k,
 \end{equation*}
 %*
where
 %*
 \begin{equation}\label{gammalp}
\widetilde{\alpha}_k=
\alpha_k+(\xi^{-1}\gamma_k-\gamma_{k+1}).\qed
 \end{equation}
 %*

Note that  $\sum_{i=1}^g \widetilde{\alpha_i} \xi^i = \sum_{i=1}^g
\alpha_i \xi^i$.

Assume that $\deg \gamma_k = 2k -1$.
 \begin{cor}
Equation \eqref{gammalp} defines a free action of the graded
additive group $\mathbb{C}^{g-1}$ with coordinates
$\gamma_2,\ldots,\gamma_g$ on the variety $V^*$. The quotient
space $V^* / \mathbb{C}^{g-1}$ is $\Gamma^*$. The vector bundle
$V^* \to \Gamma^*$ is trivial.
 \end{cor}

\prf Define the map $s:\Gamma^* \times \mathbb{C}^{g-1} \to V^*$
by the formula $ s( \xi,y, \gamma_2,\ldots,\gamma_{g-2} ) = (
\xi,t, -\gamma_2, \xi^{-1} \gamma_2-\gamma_3,\ldots, \xi^{-1}
\gamma_g )$. This is the required trivialization. \qed

Consider the case $u_1 \equiv 0$. In this case the operators
$\mathcal{U}_k$ and $\mathcal{L}$ are
 %*
 \begin{equation}\label{L0U0}
\mathcal{L}=\partial_x^2,\;\; \mathcal{U}_k=\partial_x^2
\partial_k -
\partial_{k+1}.
 \end{equation}
 %*

 \begin{lem}
Let $\alpha_1,\ldots,\alpha_g$ and $\xi$ satisfy the equation
$\left(\sum_{i=1}^g \alpha_i \xi^i\right)^2 = \xi^{-1}$. Then the
function
 %*
 \begin{equation}\label{phizero}
\Phi_0=\exp \left( \sum_{ 1 \le k \le i \le g} \alpha_{i}\, t_k\,
\xi^{i-k+1} \right)
 \end{equation}
 %*
is a common eigenfunction of the operators \eqref{L0U0} with
eigenvalues $E=\xi^{-1},\alpha_1,\ldots,\alpha_g$.
 \end{lem}

\prf The logarithmic derivative of the function  $\Phi_0$ are
given by
 %*
 \begin{equation*}
\frac{\partial_k \Phi_0}{\Phi_0}= \sum_{i=k}^g
\alpha_{i}\xi^{i-k+1}
 \end{equation*}
 %*
It is clear that $\partial_x^2 \Phi_0 = \xi^{-1} \Phi_0$ and
$\mathcal{U}_k \Phi_0 = (\xi^{-1} \partial_k - \partial_{k+1})
\Phi_0 = \alpha_k \Phi_0$. \qed

The function $\Phi_0$ can be obtained from the formula \eqref{Phi}
by rescaling. This fact will be proved in Section \ref{homcond}.

\section{Basic generating function for the solution of stationary $g$-KdV
equation} \label{basic}

Denote $\mu(\xi,\eta) = 4 \xi^{-1}+ 4 \eta^{-1}+2 \sum_{i=1}^g
\mu_{2i} \xi^i \eta^i +  \sum_{i=0}^{g-1} \mu_{2i+1}(\xi+\eta)
\xi^i \eta^i$. We have $\mu(\xi,\xi)=2 \mu(\xi)$ and
$\mu(\xi,\eta) = \mu(\eta,\xi)$, so $\mu(\xi,\eta)$ is a
polarization of $\mu(\xi)$ (see Definition \ref{defnpol}).

Consider the function
 %*
 \begin{multline*}
Q(\xi,\eta) =\mathbf{u}'(\xi) \mathbf{u}'(\eta) +
(2-\mathbf{u}(\xi)) \mathbf{u}''(\eta) +
\mathbf{u}''(\xi)(2-\mathbf{u}(\eta)) + 2(2-\mathbf{u}(\xi))
(2-\mathbf{u}(\eta)) \left(\xi^{-1}+\eta^{-1} + 2u_1\right).
 \end{multline*}
 %*

The function $Q(\xi,\eta) $ is a polarization of the function in
the right-hand side of \eqref{4mu}. Therefore $Q(\xi,\xi) =
\mu(\xi,\xi)$. Equations \eqref{polareq} and \eqref{4mu} imply
that

 %*
 \begin{equation*}
\left.\frac{ \partial \mu(\xi,\eta) }{\partial \xi}
\right|_{\xi=\eta} =\left.\frac{ \partial Q(\xi,\eta) }{\partial
\xi} \right|_{\xi=\eta}.
 \end{equation*}
 %*
Denote also
 %*
 \begin{equation*}
P(\xi) = \frac{\xi^4}{8} \left(\left.\frac{ \partial^2
\mu(\xi,\eta) }{\partial \xi \,\partial \eta} \right|_{\xi=\eta}
-\left.\frac{ \partial^2 Q(\xi,\eta) }{\partial \xi \, \partial
\eta} \right|_{\xi=\eta} \right).
 \end{equation*}
 %*

 \begin{lem}
The  function
 %*
 \begin{equation}\label{P}
P(\xi,\eta) =  \frac{ \xi^2\;\eta^2 }{4(\xi -\eta)^2} \left( 2
\mu(\xi,\eta) -Q(\xi,\eta)  \right).
 \end{equation}
 %*
is a polarization of $P(\xi)$.
 \end{lem}
\prf It is obvious that $P(\xi,\eta)$ is symmetric. Direct
calculations show that
 %*
 \begin{multline*}
P(\xi)=\frac{\xi^4}{8} \left(2 \left.\frac{\partial^2
\mu(\xi,\eta)}{\partial \xi^2}\right|_{\xi=\eta}-2
\frac{\partial^2 \mu(\xi)}{\partial \xi^2} + \left(\frac{\partial
\mathbf{u}'(\xi)}{\partial \xi} \right)^2 - 2 \frac{\partial
\mathbf{u}''(\xi)}{\partial \xi}
\frac{\partial \mathbf{u}(\xi)}{\partial \xi} \right.\\
\left. + 4 \xi^{-2}  \frac{\partial \mathbf{u}(\xi)}{\partial \xi}
(2-\mathbf{u}(\xi))+4(\xi^{-1}+u_1)\left(\frac{\partial
\mathbf{u}(\xi)}{\partial \xi}\right)^2 \right)=
\frac{1}{2}\lim_{\xi \to \eta} P(\xi,\eta).\qed
 \end{multline*}
 %*

 \begin{cor}
$P(\xi,\eta)$ is a polynomial of degree $g$ in variables  $\xi$
and $\eta$,
 \end{cor}

Define functions  $p_{ij}$ as coefficients in the expansion
 %*
 \begin{equation}\label{pij}
P(\xi,\eta) = \sum_{i=1}^{g} \sum_{j=1}^{g} p_{ij}\xi^{i}\eta^{j}.
 \end{equation}
 %*

 \begin{lem}
$p_{1i} = p_{i1} = u_i$
 \end{lem}
\prf Lemma follows from the formula $\sum_{i=1}^g p_{1i} \xi^i =
\frac{P(\xi,\eta)}{\eta}|_{\eta \to 0} = \mathbf{u}(\xi)$. \qed

This result motivates the following definition:
 \begin{defn}
The function $P(\xi,\eta)$ is called the basic generating function
for the solution $u_1$ of the stationary KdV equation.
 \end{defn}

The coefficient at $\eta^2$ in \eqref{P} is equal to :
 %*
 \begin{equation}\label{u''}
2\sum_{i=1}^g p_{2i} \xi^i  = -3 (2-\mathbf{u}(\xi))(u_1+2
\xi^{-1}) -\mathbf{u}''(\xi) + \mu_1 \xi +12 \xi^{-1}.
 \end{equation}
 %*
Therefore,
 %*
 \begin{equation}\label{u''i}
u_i''= 3 u_i u_1 +6 u_{i+1} - 2 p_{2,i}+\mu_1 \delta_{1i}
 \end{equation}
 %*
where $\delta_{ij}$ is the Kronecker symbol.

It will be shown later (see section\ref{section_sigma}) that if
$u_1=2 \wp_{gg}$ is a solution of the stationary KdV equation from
\cite{Gordon}, then $u_i=2 \wp_{g,g-i+1}$, $u_i''=2
\wp_{ggg,g-i+1}$,
$p_{2,i} = 2 \wp_{g-1,g-i+1}$. %$\mu_1=\lambda_{2g-1}$. In this case
equation  \eqref{u''i} becomes the basic relation for
$\wp$-functions (see (4.1) in \cite{Gordon}). All the results of
\cite{Gordon} for the $\wp$-functions, derived from the basic
relation, are thus true for the arbitrary solution of the
stationary KdV.

 \begin{lem}\label{dzpxi=dxipz}
$\partial(\zeta) P(\xi,\eta) = \partial(\xi) P (\zeta,\eta)$.
 \end{lem}

\prf We have
 %*
 \begin{multline*}
\partial(\zeta) P(\xi,\eta) = \frac{\xi^2\eta^2\zeta^2}{8
(\xi-\eta)(\xi-\zeta)(\eta-\zeta)} \Bigl(  \mathbf{u}''(\xi)
\bigl( \mathbf{u}'(\eta)(2-\mathbf{u}(\zeta)) -
\mathbf{u}'(\zeta)(2-\mathbf{u}(\eta))\bigr) \Bigr.\\
- \Bigl.
\mathbf{u}''(\eta)\bigl(\mathbf{u}'(\xi)(2-\mathbf{u}(\zeta)) -
\mathbf{u}'(\zeta)(2-\mathbf{u}(\xi))\bigr) +
\mathbf{u}''(\zeta)\bigl(\mathbf{u}'(\xi)(2-\mathbf{u}(\eta)) -
\mathbf{u}'(\eta)(2-\mathbf{u}(\xi))\bigr) \Bigr)\\
+\frac{ \xi\, \eta\, \zeta^2} {4
(\xi-\zeta)(\eta-\zeta)}\mathbf{u}'(\zeta)
(2-\mathbf{u}(\xi))(2-\mathbf{u}(\eta))-\frac{ \xi\, \eta^2
\zeta}{4 (\xi-\eta)(\eta-\zeta)}\mathbf{u}'(\eta)
(2-\mathbf{u}(\xi))(2-\mathbf{u}(\zeta))\\
+ \frac{ \xi^2 \eta\, \zeta}{4
(\xi-\eta)(\xi-\zeta)}\mathbf{u}'(\xi)
(2-\mathbf{u}(\eta))(2-\mathbf{u}(\zeta))\\
= \frac{1}{2}\mathcal{B}_3(\mathbf{u}''(\xi) ,
\mathbf{u}'(\xi),(2-\mathbf{u}(\xi)) + \mathcal{B}_3(
\mathbf{u}'(\xi),2-\mathbf{u}(\xi), (2-\mathbf{u}(\xi) ) \xi^{-1}
)
 \end{multline*}
 %*
Thus, $\partial(\zeta) P(\xi,\eta)$ is symmetric as a function of
variables $\xi,\eta,\zeta$. \qed

 \begin{cor}\label{cor14}
There exists a function $\phi=\phi(t_1,\ldots,t_g)$ such that
$P(\xi,\eta) = \partial(\xi)\partial(\eta) \phi$.
 \end{cor}

 \begin{cor}
$P'(\xi,\eta)=\partial(\eta)\mathbf{u}(\xi)$.
 \end{cor}
\prf Indeed,
 %*
 \begin{equation*}
P'(\xi,\eta) = \left. \frac{\partial(\zeta) P(\xi,\eta)}{\zeta}
\;\right|_{\zeta \to 0}=\frac{ \xi\;\eta }{2(\xi -\eta)}
(\mathbf{u}'(\xi)(2-\mathbf{u}(\eta))-(2-\mathbf{u}(\xi))
\mathbf{u}'(\eta))=\partial(\eta)\mathbf{u}(\xi).\qed
 \end{equation*}
 %*

Note that it follows from Theorem \ref{deuxthm} that $\partial_x
P(\xi, \eta) =\mathcal{T}_\xi^\eta  \partial_x \mathbf{u}(\xi)$.

\section{A construction of the $w$-function}

Consider the equation
 %*
\begin{equation}\label{w1}
2 \partial_x^2 \log w = - u_1,
 \end{equation}
 %*
with the initial conditions
 %*
 \begin{equation}\label{w2}
w(0)=1,\quad \partial_k w(0)=0,\; k=1,\dots,g,
 \end{equation}
 %*
here $u_1=u_1(t_1,\ldots,t_g)$ is a solution of the stationary
$g$-KdV equation with respect to $x=t_1$.

 \begin{thm}\label{wt1}
There exists a differentiable solution $w$ of \eqref{w1},
\eqref{w2} such that the functions
 %*
 \begin{equation}\label{w3}
u_k =  -2 \partial_x \partial_k \log w, \quad k=1,\ldots,g
 \end{equation}
 %*
satisfy the hypotheses of Theorem \ref{main1}.
 \end{thm}

 \begin{defn}
The solutions of \eqref{w1} described in Theorem \ref{wt1} are
called special.
 \end{defn}

\begin{thm}\label{wt3}
Let $p_{ij}(t)$ be as in \eqref{P} and \eqref{pij}. Then there is
a unique special solution of \eqref{w1} such that
 %*
 \begin{equation}\label{tt}
2 \partial_i \partial_j \log w = - p_{ij}
 \end{equation}
 %*
for all $i,j$.
 \end{thm}

\prf The existence of a required solution of \eqref{tt} follows
from Corollary \ref{cor14}. The function $w$ is defined by
\eqref{tt} up to a factor $\exp(\lambda_0+\lambda_1
t_1+\ldots+\lambda_g t_g)$. All the constants $\lambda_i$ are
uniquely determinated by the initial conditions \eqref{w2}. \qed

This result completed the solution of Problem 1.

 \begin{defn}
The special solution \eqref{w1} described in Theorem \ref{wt3} is
called a $w$-function of the solution $u$ of the stationary
$g$-KdV equation.
 \end{defn}

The relations between the higher logarithmic derivatives of the
$w$ function  are obtained with the  the technique of generating
function. For example $\sum_{ijk} \xi^i \eta^j \zeta^k \partial_i
\partial_j \partial_k \log w = \partial(\zeta) P(\xi,\eta)$. This
function was calculated in Lemma 8.3.

 The solution $u$ is a point of the space $\mathbf{R}_g$
(see Section \ref{sectionR}). Consequently, we can consider the
$w$-function as a function $w : \mathbb{C}^g \times \mathbf{R}_g
\to \mathbb{C}$.

The rest of this section is devoted to an explicit construction of
the $w$-function starting with the given solution $u$ of the KdV
equation.

Denote $t=(x,t_2,\ldots,t_g)$ and put
 %*
 \begin{equation*}
\varphi(t)= \frac{1}{2}\int_0^x \int_0^x u(t) \,\dd x
 \end{equation*}
 %*

Then equation \eqref{w1} implies: $w(t) = \exp (a(t)-\varphi(t))$
where $a''(t) = 0$. Therefore, $a(t) = a_1(\widetilde{t})x + a_0(
\widetilde{t})$ where $\widetilde{t}=(t_2,\ldots,t_g)$. The
initial condition \eqref{w2} gives now $a_0(0)=0$, $a_1(0)=0$, and
$\partial_k a_0(0)=0$, $k=2,\ldots,g$. It follows from \eqref{w3}
that
 %*
 \begin{equation}\label{w7}
2 \partial_k a_1( \widetilde{t} ) = - u_k + 2 \partial_k \int_0^x
u(t) \dd t , \quad k=2,\ldots,g.
 \end{equation}
 %*

The set of equations  \eqref{w7} with the initial condition
$a_1(0) = 0$ has a unique solution $a_1(\widetilde{ t})$. It
follows from \eqref{tt} that $2
\partial_i \partial_j a_0( \widetilde{t}) =  2 \partial _i \partial_j \varphi(t)-
2 \partial_i \partial_j \,a_1(\widetilde{ t}) x - p_{ij}(t)$.
These equations with the initial condition $a_0(0) = 0$,
$\partial_k a_0(0)=0$, $k=1,\ldots,g$ have a unique solution $a_0(
\widetilde{t})$.

\section{Applications}

\subsection{ Kleinian $\sigma$-function}
\label{section_sigma}
 Consider hyperelliptic Kleinian functions $\sigma(t)$,
$\zeta_i(t) =
\partial_i \log \sigma(t)$, and $\wp_{ij}(t)= -2 \partial_i \partial_j \log
\sigma(t)$. The function $ 2 \wp_{gg}(t)$ is a solution of the
stationary KdV equation (see \cite{Gordon}).

 \begin{cor}
Let $z \in \mathbb{C}^g$ be a point where $\sigma(z) \ne 0$. Then
the function
 %*
 \begin{equation*}
w(t) =  \frac { \sigma(t+z) } { \sigma (z)} \exp
\left<-\zeta(z),t\right>
 \end{equation*}
 %*
is a $w$-function of the solution $2 \wp_{gg}(t+z)$.
 \end{cor}

\prf The functions $u_i=2 \wp_{g,g-i+1}$ and
$p_{ij}=2\wp_{g-i+1,g-j+1}$ satisfy equations \eqref{P},
\eqref{tt} (see \cite{1}). The corollary now follows from the
uniqueness of the $w$-function. \qed

Let $\theta_g$ be the polynomials from \cite{AdlerMoser}. The
second logarithmic derivatives of $\theta_g$ give solutions of the
higher KdV equations. As it was proved in \cite{rational_abel},
the polynomial $\theta_g$ is, up to a linear change of variables,
a rational limit $\widehat{\sigma}_g$ of the $\sigma$-function of
genus $g$. Denote $\widehat{\zeta}_i(t) =
\partial_i \log \widehat{\sigma}_g(t)$.

 \begin{cor}
Let $z \in \mathbb{C}^g$ be a point where $\widehat{\sigma}_g(z)
\ne 0$. Then the function
 %*
 \begin{equation*}
w(t) = \frac { \widehat{\sigma}_g(t+z) } { \widehat{\sigma}_g (z)}
\exp \left<-\widehat{\zeta}(z),t\right>
 \end{equation*}
 %*
is the $w$-function of the solution $u = -2(\log \theta_g)''$.
 \end{cor}

\subsection{The  homogeneity  condition} \label{homcond}

The results obtained in this section follow from the uniqueness
theorems for the $w$-functions.

 \begin{lem}\label{lem24}
Suppose that $u(x,t_2,\ldots,t_g)$ is a solution of a stationary
$g$-KdV equation with respect to $x$. Take $\kappa \in
\mathbb{C}^*$. Then the function $\widehat{u}(x,t_2,\ldots,t_g) =
\kappa^2 u(\kappa x, \kappa^3 t_2,\ldots,\kappa^{2g-1} t_g)$ is
also a solution of the stationary $g$-KdV equation. Under the
transformation $u \to \widehat{u}$ the constants $\mu_i$ and $a_i$
are transformed as $\widehat{\mu}_i=\mu_i \kappa^{2i+2}$,
$\widehat{a}_i = \kappa^{2g-2i} a_i$.
 \end{lem}

\prf Let $\{u_1=u,\,u_2,\ldots,u_g\}$ be a sequence of functions
from Theorem \ref{main1}. Then the functions
 %*
 \begin{equation}\label{uihat}
\widehat{u}_i(x,t_2,\ldots,t_g) = \kappa^{2i} u_i(\kappa x
,\kappa^3 t_2,\ldots,\kappa^{2g-1}t_g)
 \end{equation}
 %*
satisfy the hypotheses of Lemma \ref{lem2} and Lemma \ref{lem3}.
 Therefore by Theorem
\ref{main1} the function $\widehat{u}$ is a solution of the
stationary KdV equation. The values $\widehat{\mu}_i$ are
determined by \eqref{4mu}, the values $\widehat{a}_i$ are found
from \eqref{mukaij}. \qed

Thus we have an action of the group $\mathbb{C}^*$ on the space
$\mathbf{R}_g$. It is obvious that under this action the initial
values $c_j=u^{(j)}(0)$ are transformed as $\widehat{c}_j =
\kappa^{2j+1} c_j$.

Denote by $\widehat{w}$ the $w$-function of the solution
$\widehat{u}$.

 \begin{lem}\label{lem25}
The $w$-functions $w$ and  $\widehat{w}$ of the solutions  $u$ and
$\widehat{u}$ are related as follows:
 %*
 \begin{equation*}
\widehat{w}(t_1,\ldots,t_g) = w(\kappa x,\kappa^3
t_2,\ldots,\kappa^{2g-1}t_g)
 \end{equation*}
 %*
 \end{lem}

\prf The functions $u_i$, $\widehat{u}_i$ are related by equation
\eqref{uihat}. It follows from \eqref{P}, \eqref{tt} that $-2
\partial_i \partial_j \widehat{w}(t_1,\ldots,t_g)= \widehat{p}_{ij}=
\kappa^{2 i+2 j-2} p_{ij} = -2 \partial_i \partial_j w(\kappa
x,\kappa^3 t_2,\ldots,\kappa^{2g-1}t_g)$. Since the $w$-function
is unique, this completes the proof.
 \qed

Consider now the $w$-function as a function on the space
$\mathbb{C}^g \times \mathbf{R}_g$.

 \begin{thm}
The function $w$ satisfies the homogeneity condition:
 %*
 \begin{multline*}
w(t_1,\ldots,t_g,a_0,\ldots,a_{g-2},c_0,\ldots,c_{2g}) = w(\kappa
t_1,\ldots,\kappa^{2g-1} t_g,\kappa^{-2g} a_0,\ldots,\kappa^{-4}
a_{g-2},\kappa^{-1} c_0,\ldots, \kappa ^{-2g-1} c_{2g})
 \end{multline*}
 %*
 \end{thm}

\prf The theorem follows directly from Lemma \ref{lem24} and Lemma
\ref{lem25}. \qed

Consider the function $\widehat{\mathbf{u}}(\xi) = \sum_{i=1}^g
\widehat{u}_i \xi_i$. It follows from \eqref{uihat} that
 %*
 \begin{equation}\label{hatbfu}
\widehat{\mathbf{u}}(x,t_2,\ldots,t_g;  \xi) =\mathbf{u}(\kappa x
,\kappa^3 t_2,\ldots,\kappa^{2g-1}t_g; \kappa^{-2} \xi).
 \end{equation}
 %*

 \begin{thm}
Let $\Phi(t_1,\ldots,t_g;\xi^{-1},\alpha_1,\ldots,\alpha_g)$ be a
common eigenfunction of the operators $\mathcal{L}$,
$\mathcal{U}_1,\ldots,\mathcal{U}_g$ with the eigenvalues
$E=\xi^{-1},\alpha_1,\ldots,\alpha_g$ (see Section
\ref{sectionF}). Then the function
 %*
 \begin{equation}\label{Phihat}
\widetilde{\Phi}(t_1,\ldots,t_g;\xi^{-1},\alpha_1,\ldots,\alpha_g;\kappa)
=\Phi(t_1 \kappa,\ldots,t_g \kappa^{2g-1};\xi^{-1} \kappa
^{-2},\alpha_1 \kappa^{-3},\ldots,\alpha_g \kappa^{-2g-1})
 \end{equation}
 %*
is regular as a function of $\kappa$ in the vicinity of the
origin, and
 %*
 \begin{equation*}
\widetilde{\Phi} = \exp\left( \sum_{1 \le i \le j \le g}
\alpha_{j}\,  \xi^{j-i+1}\, t_i\right) + O(\kappa).
 \end{equation*}
 %*
 \end{thm}

\prf Denote $\widehat{t} = (t_1 \kappa,\ldots,t_g \kappa^{2g-1})$.
Using \eqref{phiexpF}, \eqref{ditildeF}, one gets
 %*
 \begin{multline*}
\frac{\partial_i \widetilde{\Phi}}{\widetilde{\Phi}}=
\kappa^{2i-1}(\partial_i \widetilde{F})(\widehat{t};\xi^{-1}
\kappa
^{-2},\alpha_1 \kappa^{-3},\ldots,\alpha_g \kappa^{-2g-1})\\
=  \frac{ \xi^{1-i}}{4(2-u_1(\widehat{t})\xi\kappa^2-\ldots )}
\left( 2 \sum_{j \ge i} (4 \alpha_j\xi^j-u_j'(\widehat{t})\xi^j
\kappa^{2j+1})- \right.\\
\left. -\sum_{j<i} u_j(\widehat{t}) \xi^j \kappa^{2j} \sum_{j \ge
i} (4 \alpha_j\xi^j-u_j'(\widehat{t})\xi^j \kappa^{2j+1})+ \sum_{j
\ge i} u_j(\widehat{t}) \xi^j \kappa^{2j} \sum_{j < i} (4
\alpha_j\xi^j-u_j'(\widehat{t})\xi^j \kappa^{2j+1}) \right).
 \end{multline*}
 %*
Note that $\widehat{t} \to (0,\ldots,0)$ and $u_i(\widehat{t})\to
u_i(0)$ as $\kappa \to 0$. Therefore we obtain
 %*
 \begin{equation*}
\frac{\partial_i \widetilde{\Phi}}{\widetilde{\Phi}} = \sum_{j \ge
i} \alpha_j\xi^{j-i+1} + O(\kappa). \qed
 \end{equation*}
 %*

This allows to obtain a deformation of the function $\Phi$. We see
that $\widetilde{\Phi}$ tends to the function $\Phi_0$ of
\eqref{phizero} as $\kappa \to 0$.

Consider the space $L = \mathbb{C}^g \times \mathbb{C}^* \times
\mathbb{C}^g $ with coordinates $ (t_1,\ldots,t_g; \xi,
\alpha_1,\ldots,\alpha_g)$. Consider also an action of the group
$\mathbb{C}^*$ on the space $L$ given by a formula
$\kappa(t_1,\ldots,t_g;\xi,\alpha_1,\ldots,\alpha_g) =(t_1
\kappa,\, t_2 \kappa^3,\ldots,t_g \kappa^{2g-1}$; $\xi \kappa
^{2},\alpha_1 \kappa^{-3},\ldots,\alpha_g \kappa^{-2g-1})$. This
defines a projection $p\;:\; L \to M$ where $M = L/\mathbb{C}^*$.
Take some small $\varepsilon>0$ and denote $L_\varepsilon = \{
(t_1,\ldots,t_g; \xi,  \alpha_1,\ldots,\alpha_g ) \in  L : |\xi|
\ge \varepsilon\}$ and $\partial L_\varepsilon = \{
(t_1,\ldots,t_g ; \xi, \alpha_1,\ldots,\alpha_g) \in L \;:\; |\xi|
= \varepsilon \}$. Glue the boundary $\partial L_\varepsilon$ to
the space $M$ using the projection $p$ to obtains the space
$Z_\varepsilon = L_\varepsilon \cup_p M$.

Let $\varepsilon_2 < \varepsilon$. Then there is a map
$L_\varepsilon \to L_{\varepsilon_2}$ defined by the formula
$(t_1,\ldots,t_g;\xi, \alpha_1,\ldots,\alpha_g) \to
\kappa(t_1,\ldots,t_g;\xi, \alpha_1,\ldots,\alpha_g)$ where $
\kappa =\varepsilon^{-1/2}\varepsilon_2^{1/2}$. This map sends the
boundary $\partial L_\varepsilon$ to the boundary $\partial
L_{\varepsilon_2}$, so it can be lifted to a map
$Z_\varepsilon=L_\varepsilon \cup M \to
Z_{\varepsilon_2}=L_{\varepsilon_2} \cup M$. Denote
$Z=\lim_{\varepsilon \to 0} Z_\varepsilon$ and recall that
$\mathbb{C}^* \times V^*  \subset L$.

Consider the embedding $\mathbb{C}^g \times V^{*} \to Z $. This
embedding covers the embedding  $\Gamma^* \to  \Gamma$.
Approaching the limit point in $\Gamma$ corresponds to $\xi \to 0$
is the space $M \subset Z$. So we get the following result:

 \begin{thm}
On the space $Z$  there is a function $\widehat{\Phi}$, such that
$\widehat{\Phi}|_{L} = \Phi$ and $\widehat{\Phi} |_{M} = \Phi_0$.
 \end{thm}

If $\xi \to 0$ then for the restriction $\left.\Phi\right|_\Gamma
= \Phi(t_1,\ldots,t_g, \xi^{-1}, 0,0,\ldots, \alpha_g) $ one has
$\Phi \sim \exp \left( \sum_{ 1\le j  \le g} \alpha_{g}\,
\xi^{g-j+1} \, t_j \right)$.  Take a local parameter $ k =\alpha_g
\xi^g$. It follows now from the equation $ (a_g \xi^g)^2 =
\mu(\xi) = 4 \xi^{-1} + O(\xi)$ that $\Phi  \sim \exp (
\sum_{j=1}^g k^{2j-1} t_j)$.

So, the restriction $\left.\Phi\right|_\Gamma$ has the same
analytic properties as the  Baker--Akhiezer function
(\cite{krich}) of the solution $u$. By the uniqueness of the
Baker--Akhiezer function we conclude that
$\left.\Phi\right|_\Gamma$ coincides with the Baker--Akhiezer
function.

\section{Examples} %. The equations for the sequence $\{u_k\}$}
In this section we demonstrate the key constructions of the paper
in the cases $g=1$ and $g=2$.

\subsection{$g=1$}
We start with a solution $u$ of the classical stationary KdV
equation $u'''-6 u u' = 0$. Suppose that $x=0$ is a regular point
of the function $u$. Then the solution $u$ with the given values
$c_0=u(0)$, $c_1=u'(0)$, $c_2=u''(0)$ is unique in a neighbourhood
of the point $x=0$.

%Then $\mathcal{U}_1=\mathcal{A}_1=A_1$.
The key equation \eqref{4mu} becomes
 %*
 \begin{equation*}
4( 4\xi^{-1}+\mu_1 \xi + \mu_2 \xi^2) = (u')^2 \xi^2 + 2 u''\xi (
2-u \xi) + 4(\xi^{-1}+u)(2-u \xi)^2,
 \end{equation*}
 %*

hence $\mu_1 = u''-3 u^2$, $\mu_2 =1/4( (u')^2 -2 u'' u +4 u^3)$.
It is easy to see that $\mu_1'=0$ and $\mu_2'=0$. Therefore
$\mu_1$, $\mu_2$ are constants and so $\mu_1 = c_2-3 c_0^2$,
$\mu_2 =1/4( c_1^2 -2 c_2 c_0 +4 c_0^3)$.

The equation of the hyperelliptic curve is
 %*
 \begin{equation*}
4(4 \xi^{-1} + \mu_1 \xi + \mu_2 \xi^2) = y^2.
 \end{equation*}
 %*

The $w$-function is $w = \exp(-\phi(x))$ where $\phi(x) =
\frac{1}{2} \int_{0}^x \int_0^x u(x) \dd x$.

The birational equivalence $\upsilon : \mathbf{R}_1 \to
\mathcal{U}_1$ is given by the formula
$\upsilon(u_1)=\left(\Gamma, (\xi,y)\right)$, where $ \xi =
2/c_0$, $y= 2 c_1/c_0$.

Let $\Phi$ be a common eigenfunction of the operators
$\mathcal{L}$ and $\mathcal{U}_1=A_1$ with the eigenvalues
$E=\xi^{-1}$ and $\alpha$, respectively. Then the logarithmic
derivative of the function $\Phi$ is
 %*
 \begin{equation*}
\frac{\Phi'}{\Phi}=\frac{ 4 \,\alpha \,\xi - u'(x)\xi}{2( 2 -u(x)
\xi)}
 \end{equation*}
 %*
where $ 4 \alpha \xi= 2 \sqrt{4 \xi^{-1} + \mu_1 \xi + \mu_2
\xi^2}$. Therefore $\partial_x \log \Phi\sim \xi^{-1/2}$ as $\xi
\to 0$. Let $z \in \mathbb{C}$ be such that $u(z)=2 \xi^{-1}$ and
$u'(z)=4 \alpha$. Then
 %*
 \begin{equation*}
\widetilde{F}'(x;\xi^{-1},\alpha)=\widetilde{F}'(x;z)=
\frac{1}{2}\frac{u'(z)-u'(x)}{u(z)-u(x)}.
 \end{equation*}
 %*

\subsection{$g=2$}
It follows from the equation $[\mathcal{L},A_2+a_0 A_0]=0$ that
$u^{(5)}-10 u u^{(3)} -20 u'' u' +30\,u^2 u' +16  a_0 u' =0$.
Equation \eqref{mukaij} implies that $\mu_1 = 8 a_0$. We have
$\mathbf{u}(\xi) = u \xi + u_2 \xi^2$. The equation \eqref{4mu}
gives
 %*
 \begin{multline*}
4(4 \xi^{-1} + \mu_1 \xi + \mu_2 \xi^2 + \mu_3 \xi^3 + \mu_4
\xi^4)\\ = \xi ( - 12 u^2 - 16 u_2 + 4 u '') + \xi^2 ( 4 u^3 -8 u
u_2 + (u')^2 - 2 u u'' +
4 u_2 '') \\
+ \xi^3 ( 8 u^2 u_2 + 4 u_2 ^2 + 2 u' u_2' - 2 u_2 u''  - 2 u
u_2'')+ \xi^4 ( 4 u u_2^2 + (u_2')^2 - 2 u_2 u_2'').
 \end{multline*}
 %*
Therefore, $u_2 = \frac{1}{4} ( u''-3 u^2 - 8 a_0)$. Now we can
describe $\mu_2$, $\mu_3$ and $\mu_4$ as constants in the
following ordinary differential equation for $u$:
$\mu_2=\frac{1}{4} \left( 4 u^{(4)} - 10 u u'' -5 (u')^2+10 u^3
+16 a_0 u\right)$, $\mu_3 = \frac{1}{16}( 2 u' u'''-2 u u^{(4)} -2
(u'')^2 -15 u^4 + 8 u^2 u'' - 16 u^2 a_0 +12 u^2 u' +64 a_0 )$,
$\mu_4= \frac{1}{64}( (u''')^2 +16 (u'')^2 u -2 u'' u^{(4)} + 12
(u')^2 u'' + 6 u^2 u^{(4)}+ 32 u^5 - 30 u'' u^3 -12 u''' u' u -
160  a_0 u u''  + 132 a_0 u^3 + 16 a_0 u^{(4)} - 96 a_0 (u')^2 +
256 a_0^2 u )$. At last, the birational equivalence $\upsilon:
\mathbf{R}_2 \to \mathcal{U}_2$ is given by the formula
$\upsilon(u_1) =\upsilon(c_0,\ldots,c_4,a_0) = ( \Gamma,
[(\xi_1,y_1), (\xi_2,y_2)])$. Here for the construction of
$\Gamma$ the coefficients $\mu_1,\mu_2,\mu_3,\mu_4$ are used,
obtained from the formula above by substitution $c_k$ for
$u^{(k)}$. The pairs $(\xi_i,y_i)$ are the following ones:
$\xi_1,\xi_2$ are the roots of the equation $2-c_0 \xi -
\frac{1}{4}(c_2-3 c_0^2 - 8 a_0) \xi^2 =0$ and $y_1, y_2$ are
defined by the formula $y_{i} = c_1 \xi_{i} + \frac{1}{4}(c_3 - 6
c_1 c_0)$.

Acknowledgments\\
 The authors are grateful to S.P. Novikov for his attention to this
 work and stimulating discussion. We are grateful to P.G.Grinevich,
 I.M. Krichever and D.V.Leikin for valuable discussions and to
 Yu.M. Burman for  important  comments.

\end{document}